\documentclass
{elsart}

\usepackage[all]{xy}
\usepackage{graphics,graphicx}
\usepackage{epsfig}
\usepackage{amsfonts,amsmath,amssymb}

\newcommand{\be}{\begin{eqnarray*}}
\newcommand{\ee}{\end{eqnarray*}}
\newcommand{\ibe}{\begin{eqnarray}}
\newcommand{\iee}{\end{eqnarray}}

\newcommand{\pa}{\partial}
\newcommand{\f}{\frac}

\begin{document}

\begin{frontmatter}
\date{31 December 2014}

\title{Higher order operator splitting Fourier spectral methods for the Allen--Cahn equation}
\author[first]{Jaemin Shin}, \author[first]{Hyun Geun Lee}, and \author[second,cor]{June-Yub Lee}
\address[first]{Institute of Mathematical Sciences, Ewha Womans University, Seoul 120-750, Korea}
\address[second]{Department of Mathematics, Ewha Womans University, Seoul 120-750, Korea}
\ead{jyllee@ewha.ac.kr} \corauth[cor]{Corresponding author}

\begin{abstract}
The Allen--Cahn equation is solved numerically by operator splitting
Fourier spectral methods. The basic idea of the operator splitting method
is to decompose the original problem into sub-equations and
compose the approximate solution of the original equation
using the solutions of the subproblems.
Unlike the first and the second order methods, each of the heat
and the free-energy evolution operators has at least one
backward evaluation in higher order methods.
We investigate the effect of negative time steps on
a general form of third order schemes and
suggest three third order methods for better stability and accuracy.
Two fourth order methods are also presented.
The traveling wave solution and a spinodal decomposition problem
are used to demonstrate numerical properties and
the order of convergence of the proposed methods.
\end{abstract}


\end{frontmatter}

\section{Introduction} \label{intro}

The Allen--Cahn (AC) equation was originally introduced as a
phenomenological model for antiphase domain coarsening in a binary
alloy \cite{AC}:
\begin{equation}
  \frac{\partial \phi({\bf x},t)}{\partial t}
  = -\frac{F'(\phi({\bf x},t))}{\epsilon^2}
  + \Delta \phi({\bf x},t),
    \quad {\bf x} \in \Omega, ~ 0 < t \le T, \label{AC_eq}
\end{equation}
where $\Omega$ is a domain in $\mathbb R^d ~(d=1,2,3)$. The quantity
$\phi({\bf x},t)$ is defined as the difference between the
concentrations of two components in a mixture, for example,
$\phi({\bf x},t) = (m_\alpha - m_\beta)/(m_\alpha + m_\beta)$
where $m_\alpha$ and $m_\beta$ are the masses of phases $\alpha$ and $\beta$.
The function $F(\phi)=0.25(\phi^2-1)^2$ is the Helmholtz free-energy
density for $\phi$, which has a double-well form, and $\epsilon > 0$
is the gradient energy coefficient. The system is completed by
taking an initial condition $\phi({\bf x},0)=\phi^0({\bf x})$ and a
homogeneous Neumann boundary condition $\nabla \phi \cdot {\bf n}
=0$, where $\bf n$ is normal to $\partial \Omega$.

The AC equation and its various modified forms have been applied in
addressing a range of problems, such as phase transitions \cite{AC},
image analysis \cite{BCM,DB}, motion by mean curvature
\cite{ESS,KKR,FP}, two-phase fluid flows \cite{YFL}, and crystal
growth \cite{Koba,KR,BWB}. Therefore, many researchers have studied
numerical methods for solving the AC type equation
to improve stability and accuracy and to have a
better understanding of its dynamics.
Stable time step size of explicit schemes is severely restricted
due to the nonlinear term $F'(\phi)$ and implicit schemes suffer
from a solvability problem with large time steps.
One of considerable semi-implicit methods is
unconditionally gradient stable method proposed by Eyre \cite{Eyre},
which is first order accurate in time, and unconditionally gradient
stable means that a discrete energy non-increases from one time level to
the next regardless of the time step size.
And the authors in \cite{Yang,SY2} proposed first and
second order stabilized semi-implicit methods.

Another numerical method employed for solving the AC equation is
the operator splitting method \cite{Yang,LLJK,LL}. Operator
splitting schemes have been applied for many types of evolution
equations \cite{Str,GK,BC,LL2,LSL,ME}. The basic idea of the operator
splitting method is to decompose the original problem into
subproblems which are simpler than the original problem and then to
compose the approximate solution of the original problem by using
the exact or approximate solutions of the subproblems in a given
sequential order.
Operator splitting methods are simple to implement
and computationally efficient to achieve higher order accuracy
while semi-implicit schemes are hard to improve
the order of convergence.
The first and the second order operator splitting methods
for the AC equation is quite well-known \cite{Yang,LLJK,LL}, however,
the higher order (more than two) operator splitting method for the AC
equation is less well-known.

In this paper, we investigate higher order operator splitting schemes
and propose several higher order methods to solve the AC equation
with a Fourier spectral method. We
decompose the AC equation into heat and  free-energy evolution equations,
which have closed-form solutions in the Fourier and physical spaces,
respectively. Because the first and second operator splitting
methods have only forward time steps, the boundedness of the
solution is guaranteed regardless of the time step size \cite{LL}.
However, we could not guarantee the stability with large time step size
since each operator has at least one backward time step with
third and higher order of accuracy~\cite{GK,BC}.
Because a backward time marching affects numerical stability on both
sub-equations, we consider ways of minimizing the effect of
negative time steps and introduce a cut-off function
to limit the exponential amplification of high-frequency modes
in solving the heat evolution equation.

This paper is organized as follows. In section~\ref{osm}, we briefly
review the operator splitting methods which are studied by the authors in
\cite{GK}. In section~\ref{ges1}, we present higher order operator
splitting Fourier spectral methods for solving the AC equation. We
discuss the stability issues for backward time marching and suggest
the three third order operator splitting methods. We present
numerical experiments demonstrating numerical properties and
the order of convergence of the proposed methods in section~\ref{numresults}.
Conclusions are drawn in section~\ref{discu}.

\section{A brief review on the operator splitting method} \label{osm}
In this section, we review some of the basic properties of the operator
splitting  methods for a time evolution equation with two evolution terms
in summarizing the work by D.~Goldman and T.~Kapper~\cite{GK}.
Let ${\mathcal A}^{a\Delta t}$ be the solution
operator for the time evolution equation $\frac{\pa \phi}{\pa t} = f_A(\phi)$,
that is $({\mathcal A}^{a\Delta t} \phi)(t) := \phi(t + a\Delta t)$,
and ${\mathcal B}^{b\Delta t}$ be the solution operator for $f_B(\phi)$.
Then the operators ${\mathcal A}$ and ${\mathcal B}$ satisfy the semi-group properties.
Suppose we want to minimize the number of the operator evaluations of
${\mathcal A}^{a\Delta t}$ and ${\mathcal B}^{b\Delta t}$
in order to get a $N$-th order approximation of
the following ordinary differential equation
consists of two evolution terms,
\begin{equation} \label{evol-1st}
  \frac{\pa \phi}{\pa t} = f_A(\phi) +  f_B(\phi).
\end{equation}
It is well-known that the simplest form of the first order solution operator
for \eqref{evol-1st} is given as
\begin{equation} \label{s1}
  {\mathcal S}^{(1)} = {\mathcal B}^{\Delta t} \; {\mathcal A}^{\Delta t},
\end{equation}
that is, $({\mathcal S}^{(1)}\phi)(t)$ is a first order accurate approximation
of $\phi(t+\Delta t)$.
Here the choice of ${\mathcal A}$ and ${\mathcal B}$
(or $f_A$ and $f_B$) is arbitrary, thus without loss of generality, we may
assume that the first operator evaluated is always ${\mathcal A}^{a\Delta t}$.

We now consider a solution operator ${\mathcal S}^{(p)}$
with $2p$ (or $2p{-}1$ if $b_p=0$) evaluations of the operators
${\mathcal A}$ and ${\mathcal B}$ in the following form,
\begin{equation} \label{S_oper}
{\mathcal S}^{(p)} = {\mathcal B}^{b_p\Delta t} \; {\mathcal A}^{a_p\Delta t}
\; \cdots \; {\mathcal B}^{b_1\Delta t} \; {\mathcal A}^{a_1\Delta t},
\end{equation}
where all of $\{a_j\}_{j=1}^p$, $\{b_j\}_{j=1}^{p-1}$ are non-zeros.
The coefficients $a_1,\ldots,a_p$ and $b_1,\ldots,b_p$ in ${\mathcal S}^{(p)}$
must satisfy certain conditions to make ${\mathcal S}^{(p)}$
an $N$-th order approximation operator for \eqref{evol-1st}.
It is well-known that there exists ${\mathcal S}^{(p)}$ at least
$N$-th order accurate when $p \ge N$. (See \cite{GK} and
the references therein for the derivation of the following conditions.)
For first-order accuracy, $\{a_j\}$, $\{b_j\}$ must satisfy
\begin{equation} \label{1st_req}
\sum_{j=1}^{p} a_j = \sum_{j=1}^{p} b_j = 1.
\end{equation}
For second-order accuracy, $\{a_j\}$ and $\{b_j\}$
must satisfy \eqref{1st_req} and the conditions
\begin{equation} \label{2nd_req}
\sum_{j=2}^{p} a_j \left( \sum_{k=1}^{j-1} b_k \right)
= \sum_{j=1}^{p} b_j \left( \sum_{k=1}^{j} a_k \right) = \f{1}{2}.
\end{equation}
For third-order accuracy, $\{a_j\}$ and $\{b_j\}$
must satisfy \eqref{1st_req}, \eqref{2nd_req}, and the conditions
\begin{equation} \label{3rd_req}
\sum_{j=2}^{p} a_j \left( \sum_{k=1}^{j-1} b_k \right)^2
= \sum_{j=1}^{p} b_j \left( \sum_{k=1}^{j} a_k \right)^2 = \f{1}{3}.
\end{equation}

For a second-order scheme of the form \eqref{S_oper} with $p=2$,
${\mathcal S}^{(2)} = {\mathcal B}^{b_2\Delta t} \; {\mathcal A}^{a_2\Delta t}
\; {\mathcal B}^{b_1\Delta t}$  ${\mathcal A}^{a_1\Delta t}$,
\eqref{1st_req} and \eqref{2nd_req} give
\begin{equation} \label{2nd_sch}
a_1 + a_2 = 1, \quad b_1 + b_2 = 1, \quad a_2b_1 = \f{1}{2}.
\end{equation}
Since there are three equations for the four unknowns,
let $b_1{=}\omega \; (\not=0)$ be a free parameter,
then the solution of \eqref{2nd_sch} gives a general form of
a second order solution operator with up to 4 operator evaluations,
\begin{equation} \label{s2w}
{\mathcal S}_\omega^{(2)} = {\mathcal B}^{(1-\omega)\Delta t}
\; {\mathcal A}^{\frac{1}{2\omega}\Delta t}
\; {\mathcal B}^{\omega\Delta t} \; {\mathcal A}^{(1-\frac{1}{2\omega})\Delta t}.
\end{equation}
Note that ${\mathcal S}_\omega^{(2)} = {\mathcal A}^{\frac{\Delta t}{2}}
\; {\mathcal B}^{\Delta t} \; {\mathcal A}^{\frac{\Delta t}{2}}$
with $\omega=1$ is the simplest form (with only three evaluations)
among second order operators
since two evaluations of ${\mathcal A}$ and ${\mathcal B}$ is not enough
to make it second order accurate.

For a third-order scheme of the form
\begin{equation}\label{s3}
{\mathcal S}^{(3)} = {\mathcal B}^{b_3\Delta t} \; {\mathcal A}^{a_3\Delta t}
\; {\mathcal B}^{b_2\Delta t} \; {\mathcal A}^{a_2\Delta t}
\; {\mathcal B}^{b_1\Delta t} \; {\mathcal A}^{a_1\Delta t},
\end{equation}
\eqref{1st_req}, \eqref{2nd_req}, and \eqref{3rd_req} give
\begin{equation} \label{3rd_sch1}
a_1 + a_2 + a_3 = 1, \quad b_1 + b_2 + b_3 = 1,
\quad a_2b_1 + a_3(b_1+b_2) = \f{1}{2},
\end{equation}
\begin{equation} \label{3rd_sch2}
a_2b_1^2 + a_3(b_1+b_2)^2 = \f{1}{3},
\quad b_1a_1^2 + b_2(a_1+a_2)^2 + b_3 = \f{1}{3}.
\end{equation}
Choosing $b_3 = \omega$ to be a free parameter,
we can obtain two branches of the solution
for \eqref{3rd_sch1} and \eqref{3rd_sch2},
\begin{equation*}
b_1^\pm = \f{1-\omega}{2} \mp \f{\sqrt{D(\omega)}}{2(4\omega-1)}, \quad
a_2 = \f{4\omega-1}{2(3\omega-1)}, \quad
a_3^\pm = \f{\f{1}{2}-b_1^\pm a_2}{1-\omega},
\end{equation*}
\begin{equation*}
a_1^\pm = 1-a_2-a_3^\pm, \quad b_2^\pm=1-b_1^\pm-b_3,
\end{equation*}
where
\begin{equation*}
D(\omega) = (\omega-1)^2(4\omega-1)^2 + 12(4\omega-1) \left( \omega-\f{1}{3} \right)^2.
\end{equation*}
Note that real solutions of \eqref{3rd_sch1} and \eqref{3rd_sch2}
are only possible for $\omega > \f{1}{4}$ and $\omega \leq \omega^*$,
where $\omega^* \approx -1.217\cdots$ is the real root of $D(\omega)/(4\omega-1)=0$.

%
%
\begin{figure}[htb]
\center{\includegraphics[width=5in, height=3.9in]{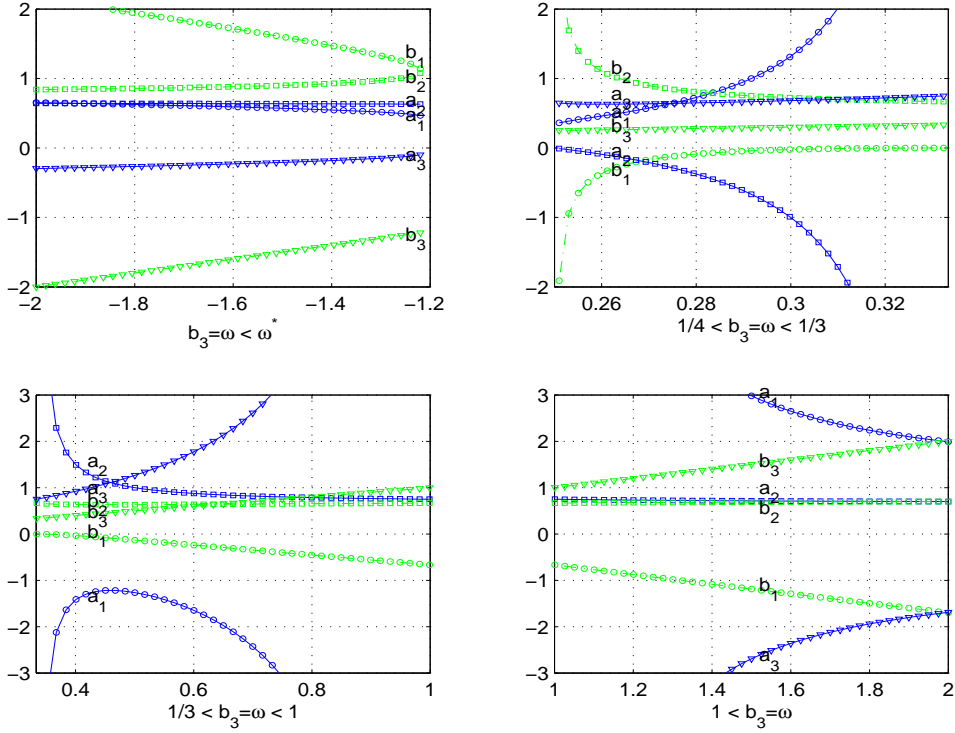}}
\caption{Positive branch solutions, $a_1^+, b_1^+, a_2^+, b_2^+, a_3^+, b_3^+$
as a function of $b_3^+=\omega$} \label{af3p}
\end{figure}
\begin{figure}[htb]
\center{\includegraphics[width=5in, height=3.9in]{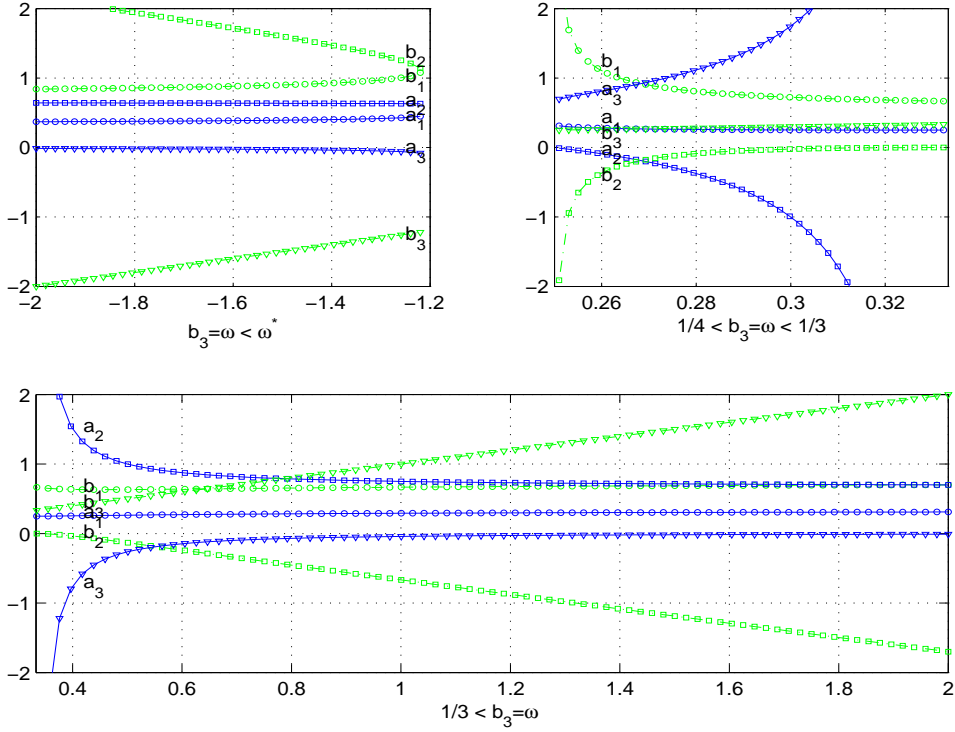}}
\caption{Negative branch solutions, $a_1^-, b_1^-, a_2^-, b_2^-, a_3^-, b_3^-$
as a function of $b_3^-=\omega$} \label{af3n}
\end{figure}

Figure~\ref{af3p} shows the positive branch solutions,
$a_1^+, b_1^+, a_2^+, b_2^+, a_3^+, b_3^+$ as a function of
$b_3^+=\omega$ for the third order operator ${\mathcal S}_{\omega^+}^{(3)}$
and Figure~\ref{af3n} shows the negative branch solutions
for ${\mathcal S}_{\omega^-}^{(3)}$.
In any case, there exists exactly one negative value among $a_1, a_2, a_3$
and also only one negative value among $b_1, b_2, b_3$.
%
%
There are three special cases when the solutions may blow up.
As $\omega \rightarrow {\f{1}{4}}^+$, ${\mathcal S}_{\omega^\pm}^{(3)}$ with
$a_2 = 0$, $b_2^\pm+b_1^\pm=\f{3}{4}$ degenerates into a second order operator,
${\mathcal B}^{\f{1}{4}\Delta t} \; {\mathcal A}^{\f{2}{3}\Delta t}
\; {\mathcal B}^{\f{3}{4}\Delta t} \; {\mathcal A}^{\f{1}{3}\Delta t}$.
As $\omega \rightarrow {\f{1}{3}}$, ${\mathcal S}_{\omega^\pm}^{(3)}$ with
$b_1^+ = 0$, $a_1^++a_2^+=\f{1}{4}$ or $b_2^-=0$, $a_2^-+a_3^-=\f{3}{4}$
degenerates into a second order operator,
${\mathcal B}^{\f{1}{3}\Delta t} \; {\mathcal A}^{\f{3}{4}\Delta t}
\; {\mathcal B}^{\f{2}{3}\Delta t} \; {\mathcal A}^{\f{1}{4}\Delta t}$.
As $\omega \rightarrow 1$, the negative branch solutions have removable
singularities and ${\mathcal S}_{\omega^-}^{(3)}$ converges to
${\mathcal B}^{\Delta t} \; {\mathcal A}^{\f{-1}{24}\Delta t}
\; {\mathcal B}^{\f{-2}{3}\Delta t} \; {\mathcal A}^{\f{3}{4}\Delta t}
\; {\mathcal B}^{\f{2}{3}\Delta t} \; {\mathcal A}^{\f{7}{24}\Delta t}$
whereas the positive branch solution does not provide a convergent operator.

We remark that a symmetric ${\mathcal S}^{(3)}$ with $b_3=0$, $a_1=a_3$, and
$b_1=b_2$ satisfying \eqref{1st_req}, \eqref{2nd_req} has only a second-order
accuracy, that is, $a_1=\f{1}{6}$, $b_1=\f{1}{2}$, and $a_2=\f{2}{3}$
does not satisfy \eqref{3rd_req}.
However, a symmetric ${\mathcal S}^{(4)}$
with $b_4=0$, $a_1=a_4$, $a_2=a_3$, and $b_1=b_3$
satisfying only \eqref{1st_req}, \eqref{2nd_req}, and \eqref{3rd_req},
\begin{equation}\label{s4w}
{\mathcal S}_U^{(4)} := {\mathcal A}^{\f{\omega}{2}\Delta t}
\; {\mathcal B}^{\omega\Delta t} \; {\mathcal A}^{\f{1-\omega}{2}\Delta t}
\; {\mathcal B}^{(1-2\omega)\Delta t} \; {\mathcal A}^{\f{1-\omega}{2}\Delta t}
\; {\mathcal B}^{\omega\Delta t} \; {\mathcal A}^{\f{\omega}{2}\Delta t}
\end{equation}
happens to be a fourth-order accuracy with
$\omega = \omega_U{=}1/(2-2^{1/3}) \approx 1.3512$,
$\frac{1-\omega}{2} \approx -0.1756$, and  $1-2\omega \approx -1.7024$.
This is the simplest form of fourth order operator with only 7 operator evaluations
and this can be derived as a symmetric combination of a second order operator
${\mathcal T}^{\Delta t} := {\mathcal A}^{\f{\Delta t}{2}}
\; {\mathcal B}^{\Delta t} \; {\mathcal A}^{\f{\Delta t}{2}}$,
\begin{equation*}
{\mathcal S}_U^{(4)} := {\mathcal T}^{\omega\Delta t} \;
{\mathcal T}^{(1-2\omega)\Delta t} \;{\mathcal T}^{\omega\Delta t},
\quad \omega = \omega_U.
\end{equation*}
Another a well-known fourth order operator splitting method~\cite{McL}
can be also derived as a symmetric combination of the second order operator
${\mathcal T}^{\Delta t}$,
\begin{eqnarray} \label{s6w}
{\mathcal S}_V^{(6)} &:= &
{\mathcal T}^{\omega\Delta t} \; {\mathcal T}^{\omega\Delta t} \;
{\mathcal T}^{(1-4\omega)\Delta t}  \;
{\mathcal T}^{\omega\Delta t} \;{\mathcal T}^{\omega\Delta t}
\\ \nonumber & = &
{\mathcal A}^{\f{\omega}{2}\Delta t} {\mathcal B}^{\omega\Delta t}
{\mathcal A}^{\omega\Delta t} {\mathcal B}^{\omega\Delta t}
{\mathcal A}^{\f{1-3\omega}{2}\Delta t} {\mathcal B}^{(1-4\omega)\Delta t}
{\mathcal A}^{\f{1-3\omega}{2}\Delta t} {\mathcal B}^{\omega\Delta t}
{\mathcal A}^{\omega\Delta t} {\mathcal B}^{\omega\Delta t}
{\mathcal A}^{\f{\omega}{2}\Delta t}
\end{eqnarray}
with $\omega = \omega_V{=}1/ (4-4^{1/3}) \approx 0.4145$.
The ${\mathcal S}_V^{(6)}$ method is computationally less
efficient (11 operator evaluations compared to minimum of 7 evaluations)
but has better stability condition
($\frac{1-3\omega}{2} \approx -0.1217$, $1-4\omega \approx -0.6580$)
than the method defined in~\eqref{s4w}.

We close this section with a remark that not just the third and
the fourth order methods mentioned above
but any operator splitting methods of third or higher order
contains at least one negative time steps for each of the operators,
${\mathcal A}$, ${\mathcal B}$.
(See \cite{GK,BC} for the proof of the general theorem.)

\section{Higher-order operator splitting Fourier spectral methods} \label{ges1}

We consider the AC equation (\ref{AC_eq}) in one-dimensional space
$\Omega=(0,L)$. Two- and three-dimensional spaces can be analogously defined.
For simplicity of notation, we sometimes abuse the notation
$\phi = \phi(t)$ referring $\phi(\cdot, t)$ and
define the ``{\em free-energy evolution operator}''
${\mathcal F}^{\Delta t}$ as follows
\begin{equation}
{\mathcal F}^{\Delta t} (\phi(t^n)) := \phi(t^n + \Delta t),
\end{equation}
where $\phi(t^n + \Delta t)$ is a solution of the first order differential
equation
\begin{equation*}
   \f{\pa \phi}{\pa t} = - \f{F'(\phi)}{\epsilon^2}
\end{equation*}
with an initial condition $\phi(t^n)$. For given $F'(\phi)=\phi^3-\phi$, we
have an analytical formula (See~\cite{Yang,LLJK,LL}) for the evolution operator
${\mathcal F}^{\Delta t}$ in the physical space
\begin{equation} \label{evol-free}
   {\mathcal F}^{\Delta t} (\phi)
   = \frac{\phi} {\sqrt {\phi^2 + (1-\phi^2) e^{-\frac{2\Delta t}{\epsilon^2}} }}.
\end{equation}

We also define the ``{\em heat evolution operator}'' ${\mathcal H}^{\Delta t}$
as follows
\begin{equation}
{\mathcal H}^{\Delta t} (\phi(t^n)) := \phi(t^n + \Delta t),
\end{equation}
where $\phi(t^n + \Delta t)$ is a solution of the first order differential
equation
\begin{equation*}
\f{\pa \phi}{\pa t} = \Delta \phi
\end{equation*}
with an initial condition $\phi(t^n)$. In this paper, we employ the discrete
cosine transform \cite{NTK} to solve the AC equation with the zero Neumann
boundary condition: for $k=0,\ldots,M{-}1$,
\begin{equation*}
\widehat{\phi}_{k}
   = \alpha_{k} \sum_{l=0}^{M{-}1}
   \phi_{l} \cos \left[ \f{\pi}{M} k
      \left( l{+}\f{1}{2} \right) \right],
\end{equation*}
where $\phi_{l} = \phi\left( \f{L}{M}\left( l{+}\f{1}{2} \right) \right)$ and
$\alpha_{0} = \sqrt{1/M}$, $\alpha_{k} = \sqrt{2/M}$ for $k \ge 1$. Then, we
have a semi-analytical formula for the evolution operator ${\mathcal H}^{\Delta t}$
in the discrete cosine space
\begin{equation}
   {\mathcal H}^{\Delta t} (\phi)
   = {\mathcal C}^{-1}\left[ \: e^{A_{k} \Delta t}
     {\mathcal C} \left[ \phi \right] \right],
\end{equation}
where $A_{k}=- \left( \f{\pi k}{L} \right)^2$ and ${\mathcal C}$ denotes the
discrete cosine transform.

For the first order operator splitting scheme ${\mathcal S}^{(1)}$ in~\eqref{s1}
and the second order scheme ${\mathcal S}_\omega^{(2)}$
in \eqref{s2w} with $0 < \omega  \le 1$,
the evaluations are all forward time marching, that is,
all of $\{a_j\}_{j=1}^p$ and $\{b_j\}_{j=1}^p$ are positive.
We can easily show that both schemes are unconditionally stable, in the sense that
$|\phi(t^n+\Delta t)| \leq 1$ if $|\phi(t^n)| \leq 1$ regardless of the time
step size. (See \cite{LL} for the proof.)
However, in the case of third or higher order, each of
operators ${\mathcal F}$, ${\mathcal H}$ has
at least one backward evaluation as mentioned in section~\ref{osm}.
For this reason, we need to investigate the stability of
the operators ${\mathcal F}^{-\Delta t}$ and ${\mathcal H}^{-\Delta t}$
especially for large $\Delta t$.

The stability issue for backward time heat equation is well-known.
Even though ${\mathcal H}^{\pm\Delta t} \; {\mathcal H}^{\mp\Delta t}$
(without noise) is always an identity operator regardless of the size of $\Delta t$,
the numerical composition of the operators (even with small error)
is far away from the identity operator when $\Delta t$ becomes large
since ${\mathcal H}^{-\Delta t}$ is exponentially big for $\Delta t \gg 1$.
The stability of the numerical composition of the free energy evolution
operators is less well-known and we want to explain
why the numerical composition of the operators
${\mathcal F}^{\pm\Delta t} \; {\mathcal F}^{\mp\Delta t}$
(even with small error) is far away from the identity operator
when $\Delta t$ becomes large using the following figure.

\begin{figure}[h]
\center{\includegraphics[width=5in]{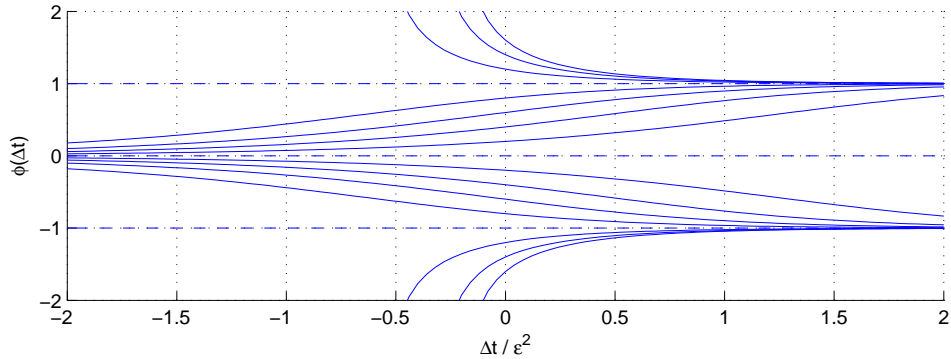}}
\caption{$\phi(\Delta t) = {\mathcal F}^{\Delta t} (\phi(0))$
with various initial values, $\phi(0) = -1.6, -1.4, \cdots, 1.6$}
\label{freeT}
\end{figure}
Figure~\ref{freeT} plots ${\mathcal F}^{\Delta t} (\phi)$ as a function of
$\Delta t$ with various initial values of $\phi$ between $-1.6$ to $1.6$.
As you can see, ${\mathcal F}^{\Delta t} (\phi)$ with $|\phi| < 1$
converges to $\pm 1$ as $\Delta t \gg 1 $, however,
the solution with $|\phi| > 1$ as a result of small perturbation
may blow up when $\Delta t \ll -1$.
Therefore, composition of two operators
${\mathcal F}^{\Delta t}$ followed by ${\mathcal F}^{-\Delta t}$
even with small evaluation error near 1 is no longer bounded
as $\Delta t$ is getting bigger.
And ${\mathcal F}^{-\Delta t} (\phi)$ with $|\phi| < 1$ converges to 0
for $\Delta t \gg 1$ thus ${\mathcal F}^{-\Delta t}$ followed by
${\mathcal F}^{\Delta t}$ for $\Delta t \gg 1$
may change the sign of result even with small perturbation near 0.
This non-linear stability effect is basically a consequence of
the fact that the solution of the free-energy evolution operator
${\mathcal F}^{\Delta t} (\phi)$
is exponentially close to 1 or 0 as $\Delta t \rightarrow \pm \infty$.

%
%
\begin{figure}[h]
\center{\includegraphics[width=5in]{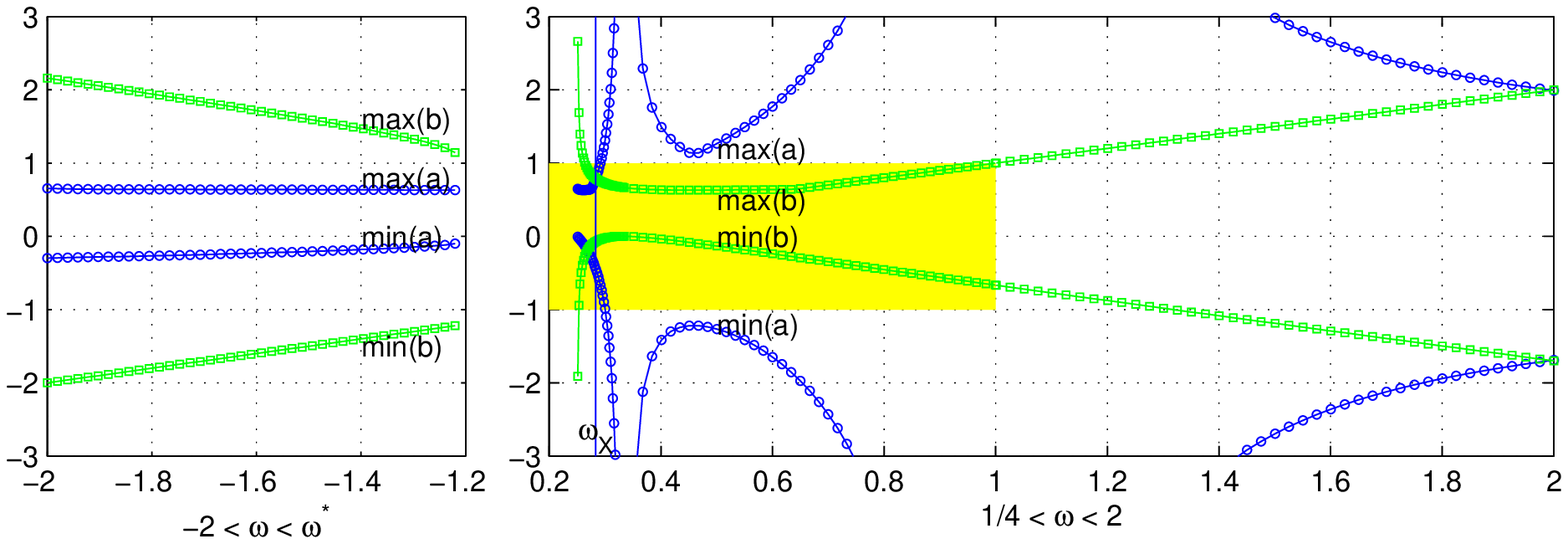}} \caption{Minimum and maximum
of $\{a_i^+\}_{i=1}^3$ and $\{b_j^+\}_{j=1}^3$ as a function of $b_3^+=\omega$.
The region where values are bounded by $[-1,1]$ is shaded in yellow.}
\label{aw3p}
\end{figure}
\begin{figure}[h]
\center{\includegraphics[width=5in]{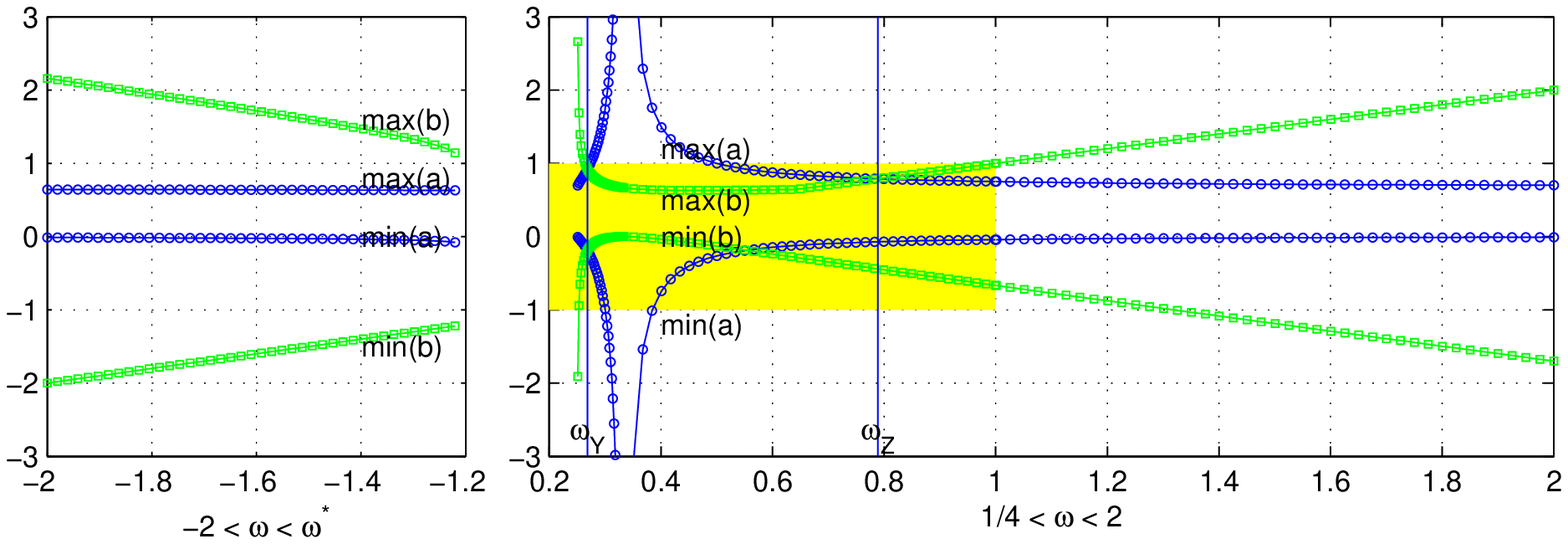}} \caption{Minimum and maximum
of $\{a_i^-\}_{i=1}^3$ and $\{b_j^-\}_{j=1}^3$ as a function of $b_3^-=\omega$.
The region where values are bounded by $[-1,1]$ is shaded in yellow.}
\label{aw3n}
\end{figure}

In order to achieve a better stability condition,
we propose third order schemes with bounded values of $\{a_j, b_j\}_{j=1}^p$.
Figures \ref{aw3p} and \ref{aw3n} show the minimum and maximum values
of $\{a_j\}_{j=1}^3$ and $\{b_j\}_{j=1}^3$ for the positive and negative
branches, respectively.
It is worth noting that $a_1^+, b_1^+, a_2^+, b_2^+, a_3^+, b_3^+$
are bounded by $[-1, 1]$ when $ 0.26376\cdots \le \omega^+ \le 0.29167\cdots$ and
$a_1^-, b_1^-, a_2^-, b_2^-,a_3^-,  b_3^-$ are bounded by $[-1, 1]$
when $ 0.26376\cdots \le \omega^- \le  0.27362\cdots$ or $ 1/2 \le \omega^- \le 1$.
Since there is exactly one negative value among $\{a_j\}_{j=1}^3$,
$\max\{a_j\}_{j=1}^3 \ge -\min\{a_j\}_{j=1}^3$
can be inferred from \eqref{1st_req} when $|a_j| \le 1$.
Similarly $\max\{b_j\}_{j=1}^3 \ge -\min\{b_j\}_{j=1}^3$
when $|b_j| \le 1$.
In the shaded regions on the figures
where values of $|a_j| , |b_j|$ are bounded by $1$,
there are three local minima of $\max\{|a_j|, |b_j|\}_{j=1}^3$
at which values are summarized on the following table.
\begin{table}[h]
\caption{Solutions for ${\mathcal S}_{\omega^\pm}^{(3)}$,
$a_1^\pm, b_1^\pm, a_2^\pm, b_2^\pm, a_3^\pm, b_3^\pm$
at the local minima of $\max\{|a_j|, |b_j|\}_{j=1}^3$.\label{table3rd}}
\begin{center}
\begin{scriptsize}
\begin{tabular} {|c|c|c|c|c|c|c|c|c|}
\hline
$\omega^\pm$ &  Condition    & $a_1$     & $b_1$      & $a_2$      & $b_2$      & $a_3$      & $b_3$      \\
\hline
$\omega_X$ & $a_1^+ = b_2^+$ & 0.78868.. & -0.07189.. & -0.44191.. &  0.78868.. &  0.65324.. &  0.28322.. \\
$\omega_Y$ & $b_1^- = a_3^-$ & 0.26833.. &  0.91966.. & -0.18799.. & -0.18799.. &  0.91966.. &  0.26833.. \\
$\omega_Z$ & $a_2^- = b_3^-$ & 0.28322.. &  0.65324.. &  0.78868.. & -0.44191.. & -0.07189.. &  0.78868.. \\
\hline
\end{tabular}
\end{scriptsize}
\end{center}
\end{table}


It is worth to note that the sets of $\{a_j^-\}_{j=1}^3$
and $\{b_j^-\}_{j=1}^3$ are same for $\omega^- = \omega_Y$.
The set $\{a_j^+\}_{j=1}^3$ for $\omega^+ = \omega_X$
is $\{b_j^-\}_{j=1}^3$ for $\omega^- = \omega_Z$ and
the set $\{b_j^+\}_{j=1}^3$ for $\omega^+ = \omega_X$
is $\{a_j^-\}_{j=1}^3$ for $\omega^- = \omega_Z$.
This symmetry gives us a freedom to choose the order of operator
evaluations and we define three third order operator splitting methods
${\mathcal S}_X^{(3)},  {\mathcal S}_Y^{(3)}, {\mathcal S}_Z^{(3)}$
for the AC equation as follows:
\begin{equation} \label{S3XYZ}
{\mathcal S}_X^{(3)},  {\mathcal S}_Y^{(3)}, {\mathcal S}_Z^{(3)}
:= {\mathcal F}^{b_3\Delta t} \; {\mathcal H}^{a_3\Delta t}
  \; {\mathcal F}^{b_2\Delta t} \; {\mathcal H}^{a_2\Delta t}
  \; {\mathcal F}^{b_1\Delta t} \; {\mathcal H}^{a_1\Delta t}
\end{equation}
where $\{a_j\}_{j=1}^3$ and $\{b_j\}_{j=1}^3$ are given in Table~\ref{table3rd}.

Another issue raised with negative time step is that
the heat evolution operator ${\mathcal H}^{a_j \Delta t},  a_j<0$  may
amplify the high frequency modes exponentially big, $e^{A_k a_j \Delta t} \gg 1$.
This situation $-A_k \Delta t = \left(\frac{\pi k}{L}\right)^2 \Delta t \gg 1$
happens when $k^2 \Delta t \gg O(1)$. On the other hand,
a physically reasonable bound for $\Delta t$ in the AC equation
is $\frac{\Delta t}{\epsilon^2} \le O(1)$, thus
the blow-up may occur only for physically too high
frequency modes, $k \gg \frac{L}{\epsilon}$.
Thus, we introduce a cut-off function to bound of
${\mathcal H}^{a_j \Delta t}$ for high frequency modes
where $-A_k \Delta t \gg 1$.
We will numerically demonstrate the effect of introducing
the cut-off function in subsection~\ref{subsec-cutoff}.

\section{Numerical experiments} \label{numresults}

In this section, we numerically demonstrate the order of
convergence of the proposed third order schemes ${\mathcal S}_X^{(3)},
{\mathcal S}_Y^{(3)},  {\mathcal S}_Z^{(3)}$ in \eqref{S3XYZ}
and the fourth order schemes ${\mathcal S}_U^{(4)}$ in \eqref{s4w}
and ${\mathcal S}_V^{(6)}$ in \eqref{s6w}.
Two examples are used for the test, one is the traveling wave
solution with analytic solution and the other is
a three-dimensional spinodal decomposition problem with random initial values.

\begin{figure}[htbp]
\centering
\includegraphics[width=2.5in]{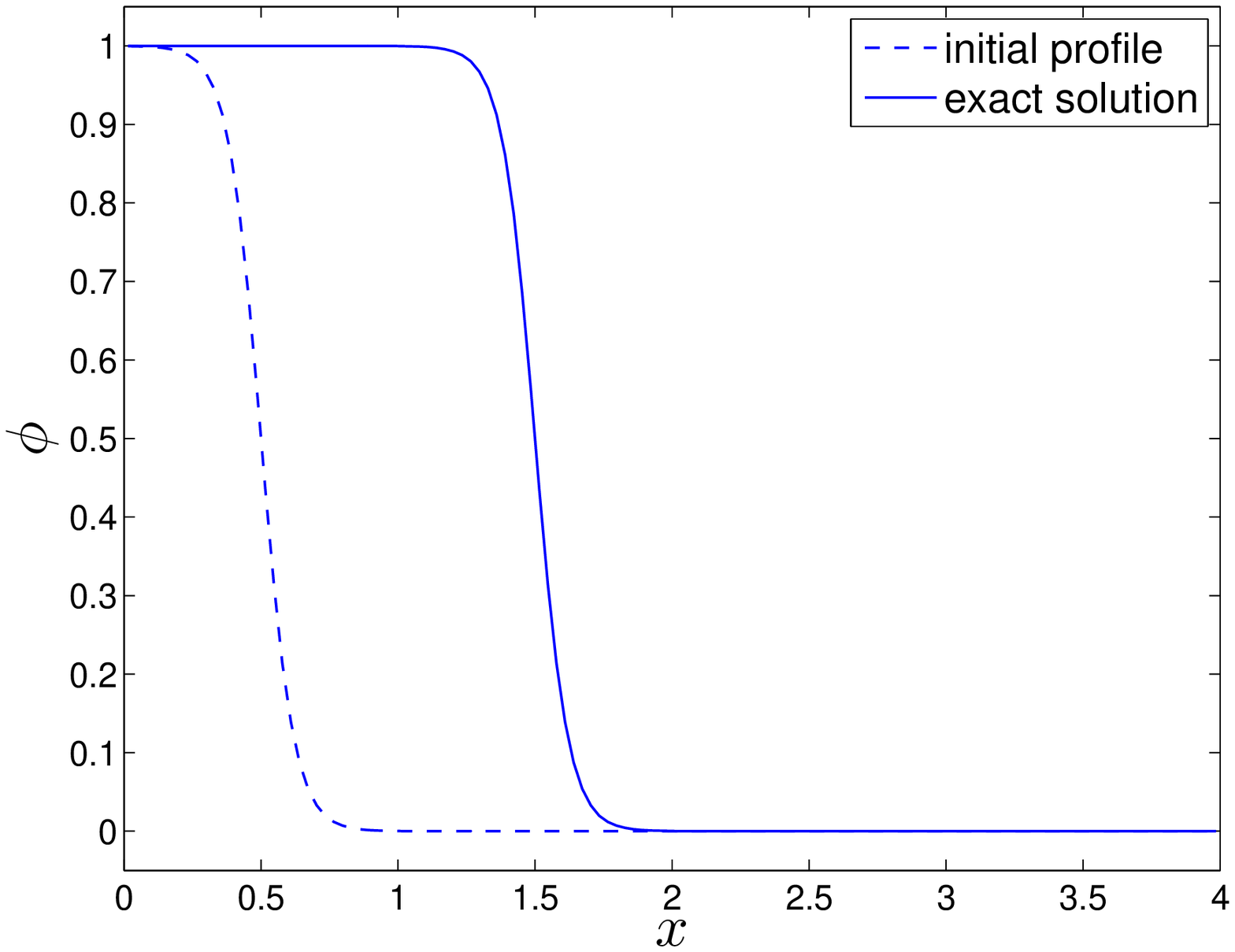}
\includegraphics[width=2.5in]{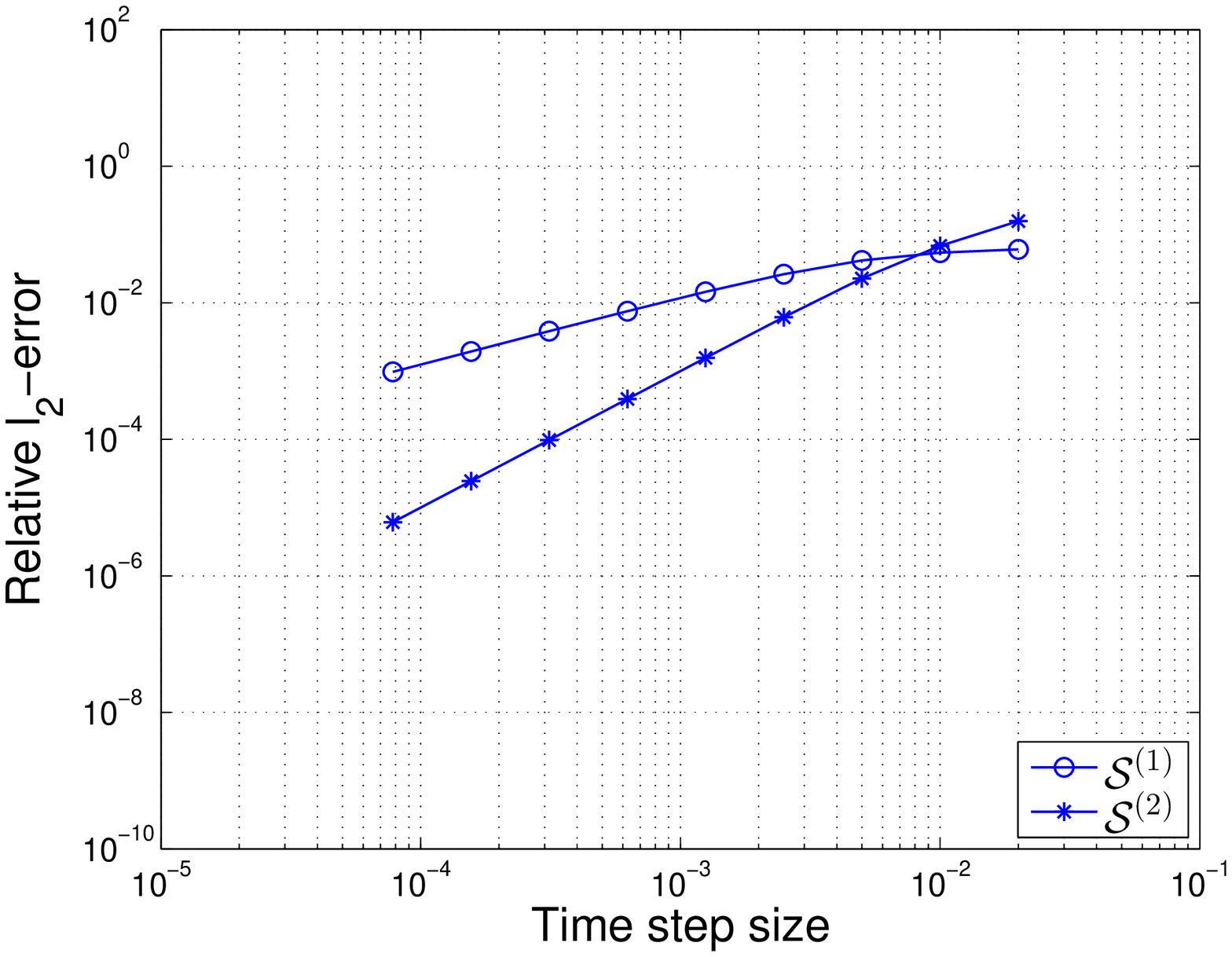}
\caption{Traveling wave solution $\phi(x,T_f)$ at $T_f=1/s$
with $\epsilon=0.03\sqrt{2}$.
And relative $l_2$ errors of $\phi(x,T_f)$ by $\mathcal S^{(1)},
\mathcal S^{(2)}$ with $h=2^{-5}$ for various time step sizes
$\Delta t$.} \label{acc_test12}
\end{figure}

One of the well-known traveling wave solutions of the Allen--Cahn equation is
\begin{equation} \label{wave1d}
 \phi(x,t) = \f{1}{2} \left(
1-\tanh\f{x-0.5-st}{2\sqrt{2}\epsilon} \right),
\end{equation}
where $s=3/(\sqrt{2}\epsilon)$ is the speed of the traveling wave.
The leftmost plot in Figure~\ref{acc_test12} shows
the initial profile $\phi(x,0)$ and the analytic solution
$\phi(x,T_f)$ at $T_f = 1/s$ with $\epsilon=0.03\sqrt{2}$.
Using this traveling wave solution,  we compare the first, second,
third, and fourth order operator splitting Fourier spectral
methods described in section 3.
The numerical solutions $\phi(x,t), \; 0 < t \le T_f$
are obtained with various time step sizes $\Delta t$
but the spatial grid size is fixed to $h=2^{-5}$ which
provides enough spatial accuracy.
The traveling wave solution with the same numerical parameters
are used in the following two subsections
to test the third and the fourth order schemes.

The rightmost plot in Figure~\ref{acc_test12} shows
the numerical error of the first order scheme $\mathcal S^{(1)}$ in \eqref{s1}
and the second order scheme $\mathcal S_{\omega=1}^{(2)}$ in \eqref{s2w}
compared to the analytic solution at $t=T_f$.
It is worth to remind that the first and the second order schemes
apply only forward time steps of ${\mathcal F}$ and ${\mathcal H}$,
thus the stability (or boundedness of the solution)
regardless of the size $\Delta t$ can be easily proven.
(See our previous paper~\cite{LL} for numerical properties of
these first and second order schemes.)

\subsection{Cut-off function and stability of the third order methods}
\label{subsec-cutoff}
As mentioned in section 3, negative time steps of ${\mathcal F}$
and ${\mathcal H}$ are unavoidable in the third or higher order operator
splitting methods. Especially a negative time step makes the heat evolution
operator exponentially big, therefore,
we introduce the following cut-off function with a tolerance $K_{tol}$
for the heat evolution operator ${\mathcal H}$,
\begin{equation} \label{cutoff}
{\mathcal H}^{a_j \Delta t} (\phi) =
{\mathcal C}^{-1}\left[ \: \min \{ e^{A_{k} a_j \Delta t}, K_{tol} \}
    \;\; {\mathcal C} \left[ \phi \right] \right].
\end{equation}
The choice of $K_{tol}$ depends on the time step size
$a_j \Delta t$ and highest frequency modes $k_{max}$
which are functions of desired computational accuracy.
Following computational examples in this subsection give
a basic guideline for the choice of $K_{tol}$.

\begin{figure}[htbp]
\centering
\includegraphics[width=4.5in]{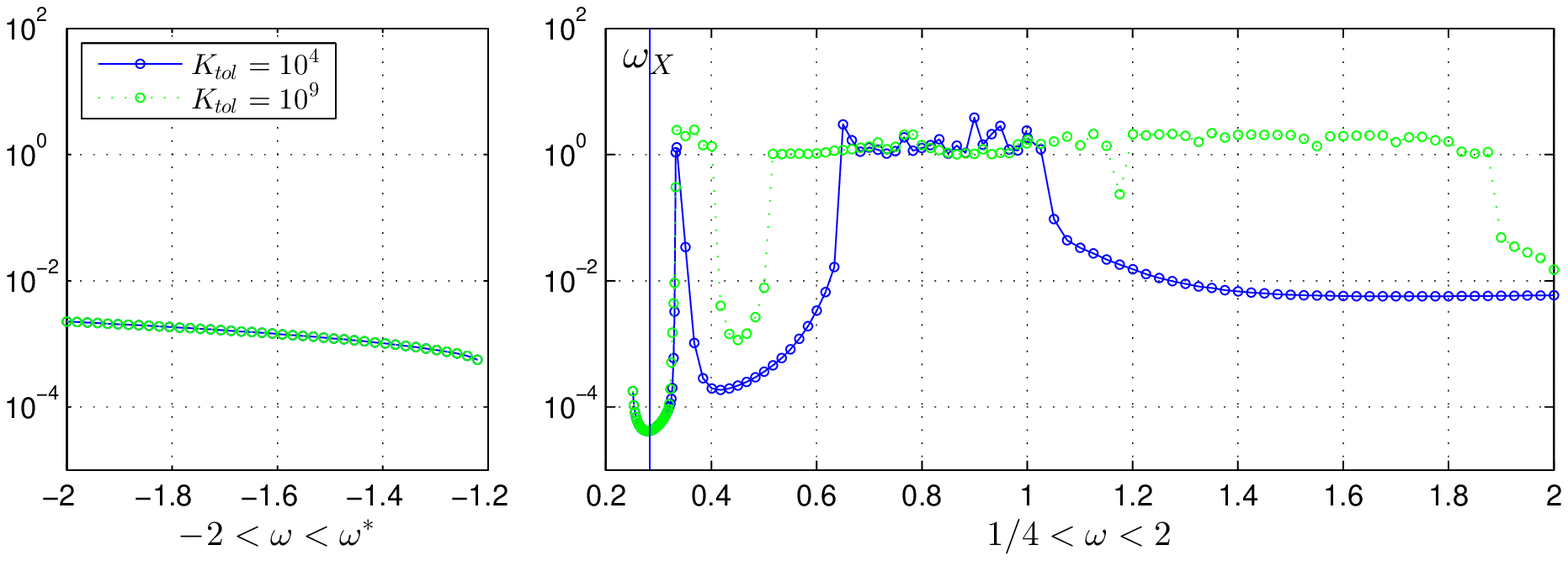} \\
(a) With positive branch solutions for ${\mathcal S}_\omega^{(3)}$ in \eqref{s3}
\includegraphics[width=4.5in]{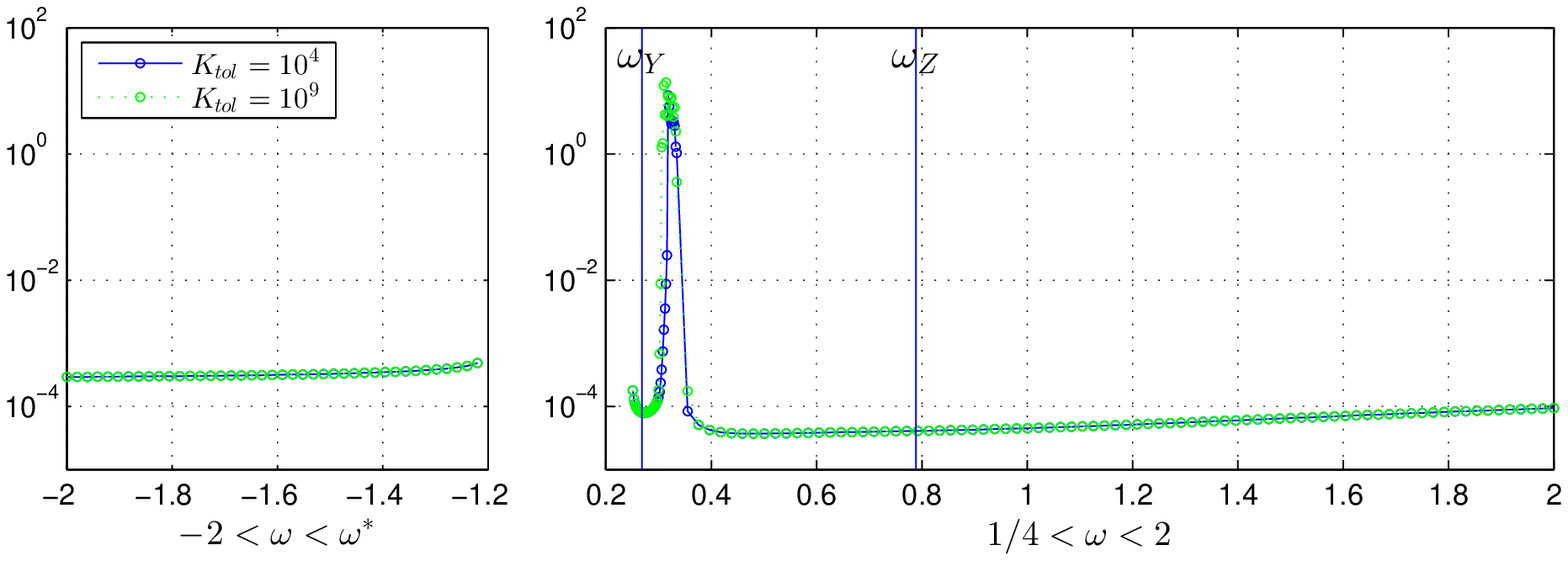} \\
(b) With negative branch solutions for ${\mathcal S}_\omega^{(3)}$ in \eqref{s3}
\caption{Relative $l_2$ errors of the traveling wave solution
$\phi(x,T_f=1/s)$ by the third order method ${\mathcal S}_\omega^{(3)}$
for various $\omega$  with $\Delta t=2^{-4}/s$,
$\epsilon=0.03\sqrt{2}$, and $h=2^{-5}$.}
\label{ome_eff}
\end{figure}

As mentioned in section \ref{osm}, we have various coefficients
$\{ a_j^\pm \}_{j=1}^3$ and $\{b_j^\pm \}_{j=1}^2$ as a function
of $b_3^\pm=\omega$.
To investigate the effect of $\omega$ in the third order method
${\mathcal S}_\omega^{(3)}$,
we consider the traveling wave problem given in~\eqref{wave1d}.
We compute relative $l_2$ errors  for various $\omega$ values
with a fixed time step $\Delta t=2^{-4}/s$ and
Figures~\ref{ome_eff} (a) and (b) show relative $l_2$ errors of the
traveling wave solution $\phi(x,T_f)$ by the third order methods
${\mathcal S}_\omega^{(3)}$
for positive and negative branches of various $\omega$, respectively.
Here we set $K_{tol}=10^4$ (blue solid line) or $10^9$ (green dashed line).

The first noticeable point in Figure~\ref{ome_eff} might be
that the error is relatively large at $\omega^\pm \to {\f{1}{4}}^+$
or $\omega^\pm \to \f{1}{3}$ where the third order operator
degenerates into a second order operator.
Also a region near $\omega^+=1$ in the positive branch case,
the computation does not provide any accuracy at all.
As $\omega^+ \to 1$, ${\mathcal S}_\omega^{(3)}$ contains
a big negative time step of the heat evolution operator
${\mathcal H}^{a_j \Delta t}$ since $\min\{a_j\} \to -\infty$.
In these cases, the choice of cut-off parameter $K_{tol}$ becomes
important, and small $K_{tol}$ is recommended when
$-\min\{a_j\} \Delta t \gg \left( \frac{L}{\pi k_{max}} \right)^2$.

For $ 0.26376\cdots \le \omega^+ \le 0.29167\cdots$
in which $\{a_j^+\}$ and $\{b_j^+\}$ are bounded by $[-1, 1]$,
especially near $\omega_X$ at which $\max \{ |a_j|, |b_j| \}$ has a local minimum,
the error is smaller than that for other $\omega$ values.
The similar phenomenon is observed the computation for the negative branch.
We choose three special values $\omega^+ = \omega_X$, $\omega^- = \omega_Y$,
and $\omega^- = \omega_Z$ for ${\mathcal S}_X^{(3)}$, ${\mathcal S}_Y^{(3)}$,
and ${\mathcal S}_Z^{(3)}$, respectively. For these cases,
all $\{a_j\}$ are bounded by $[-1,1]$ and the choice
of cut-off value $K_{tol}$ does not play an important role
in the computation.

\begin{figure}[htbp]
\begin{minipage}{0.49\linewidth}
\centering
\includegraphics[width=2.5in]{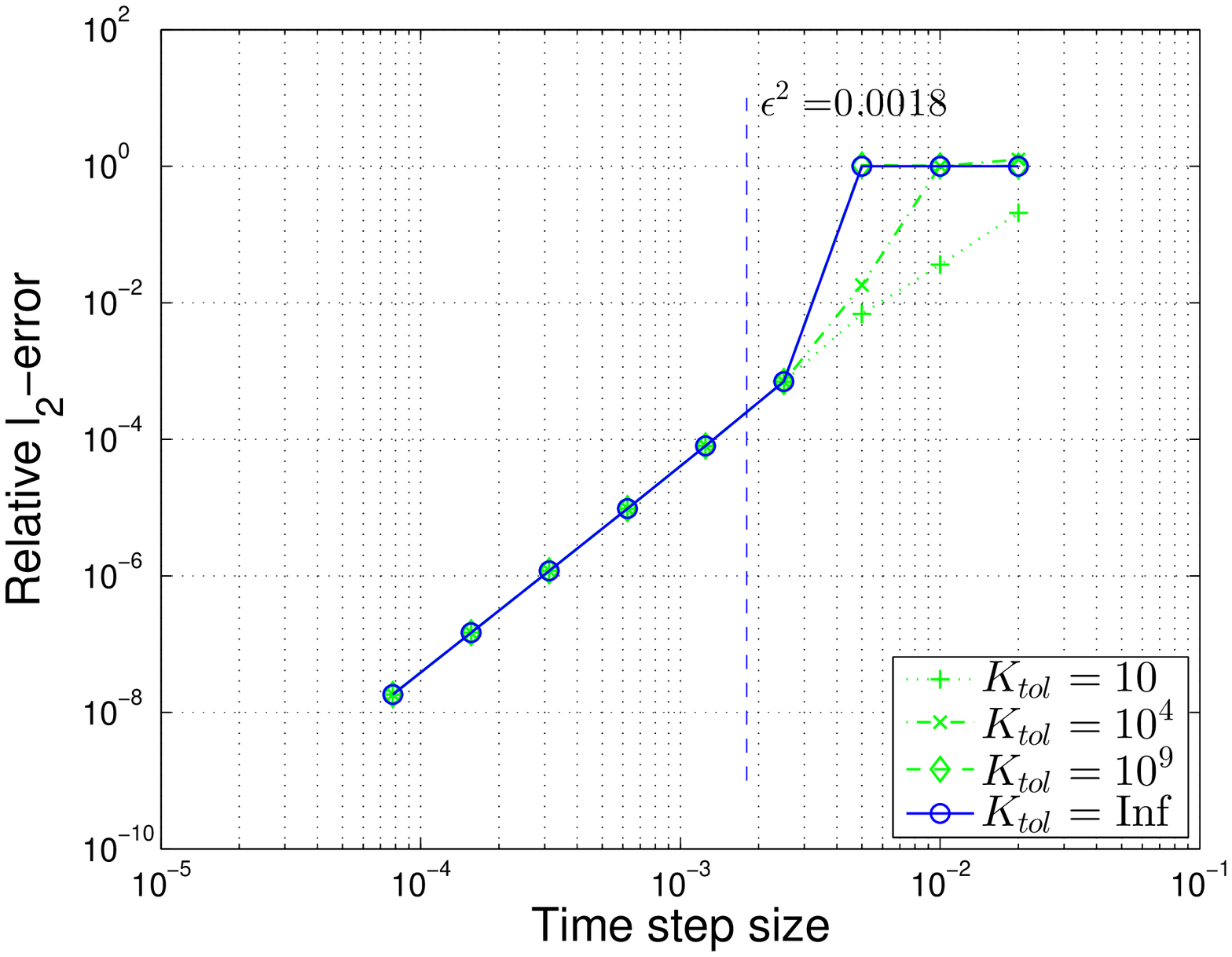} \\
(a) $M=256$
\end{minipage}
\begin{minipage}{0.5\linewidth}
\centering
\includegraphics[width=2.5in]{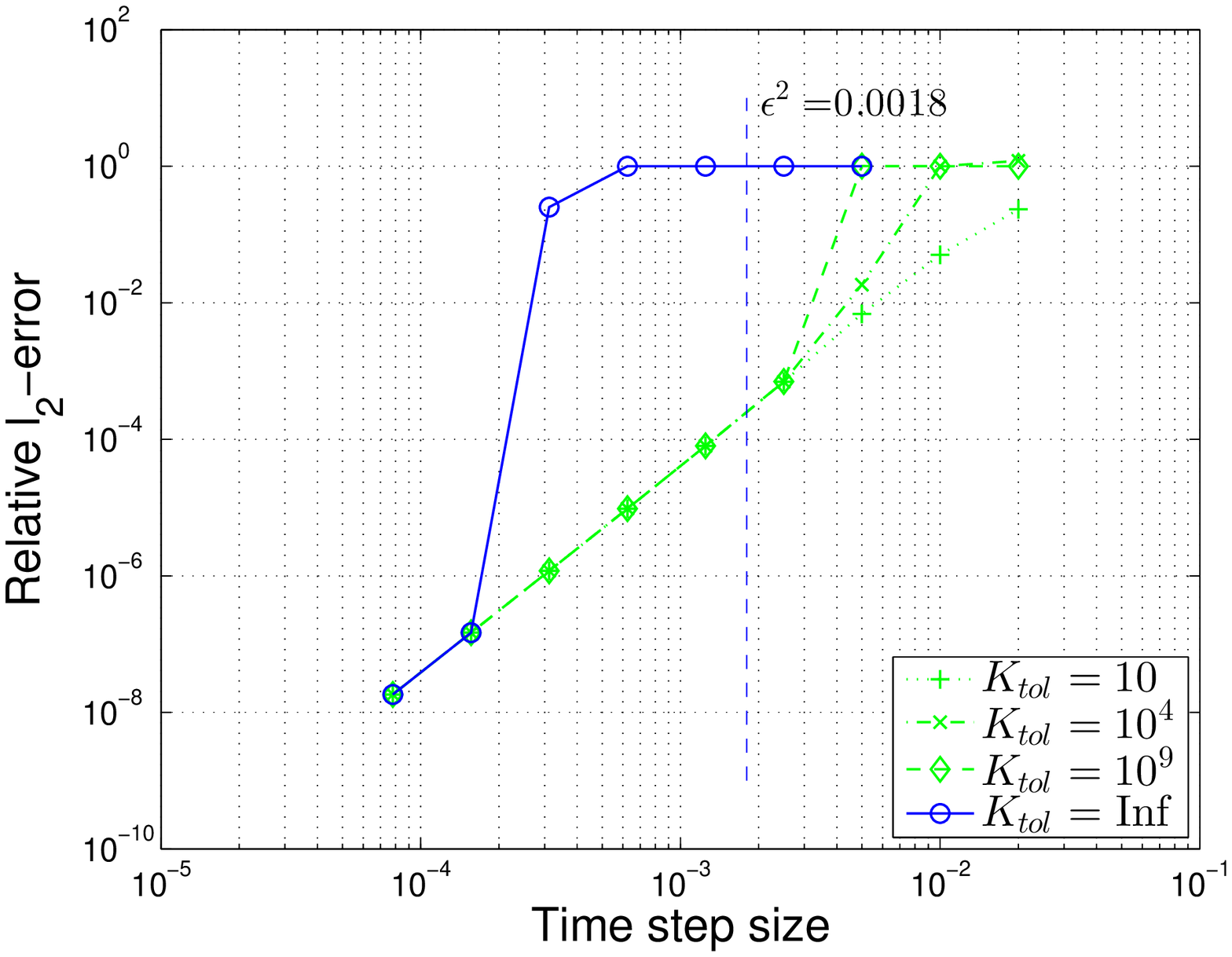} \\
(b) $M=1024$
\end{minipage}
\caption{Relative $l_2$ errors of the traveling wave solution $\phi(x,T_f=1/s)$
by ${\mathcal S}_Y$ with $\epsilon=0.03\sqrt{2}$.}
\label{cutoff_nx}
\end{figure}

We now investigate the effect of highest frequency $k_{max}$ to $K_{tol}$.
Plots in Figure~\ref{cutoff_nx} show relative $l_2$ errors of
the traveling wave solution $\phi(x,T_f)$ by the third order method
${\mathcal S}_Y^{(3)}$ with different spatial grid sizes
$h = \frac{L}{M} = \frac{4}{256} = 2^{-6}$ or $h=\frac{4}{1024} = 2^{-8}$.
If a cut-off function is not used (labeled as $K_{tol}=\mbox{Inf}$),
the computation provides no accuracy for relatively large time step.
The computation may even stop as two biggest $\Delta t$ cases for $M=1024$
and the cases happen more often as $k_{max} = M$ becomes large.
If $\Delta t$ is larger than $\epsilon^2$,  $K_{tol}$ must be properly chosen
in order to valence the accuracy loss while avoiding blow-up.
However, the choice of $K_{tol}$ makes no significant difference
of the solution when $\Delta t \le \epsilon^2$ (which is physically valid limit
for the AC equation) since the high frequency modes ${\widehat\phi}_k$
with $k \gg \frac{L}{\epsilon}$ is negligible for the physically meaningful
solution of the AC equation. So the simplest rule of thumb might be
setting $K_{tol}$ around the desired accuracy of the computation.

\subsection{Convergence of the third and the fourth order methods}
\begin{figure}[htbp]
\begin{minipage}{0.49\linewidth}
\centering
\includegraphics[width=2.5in]{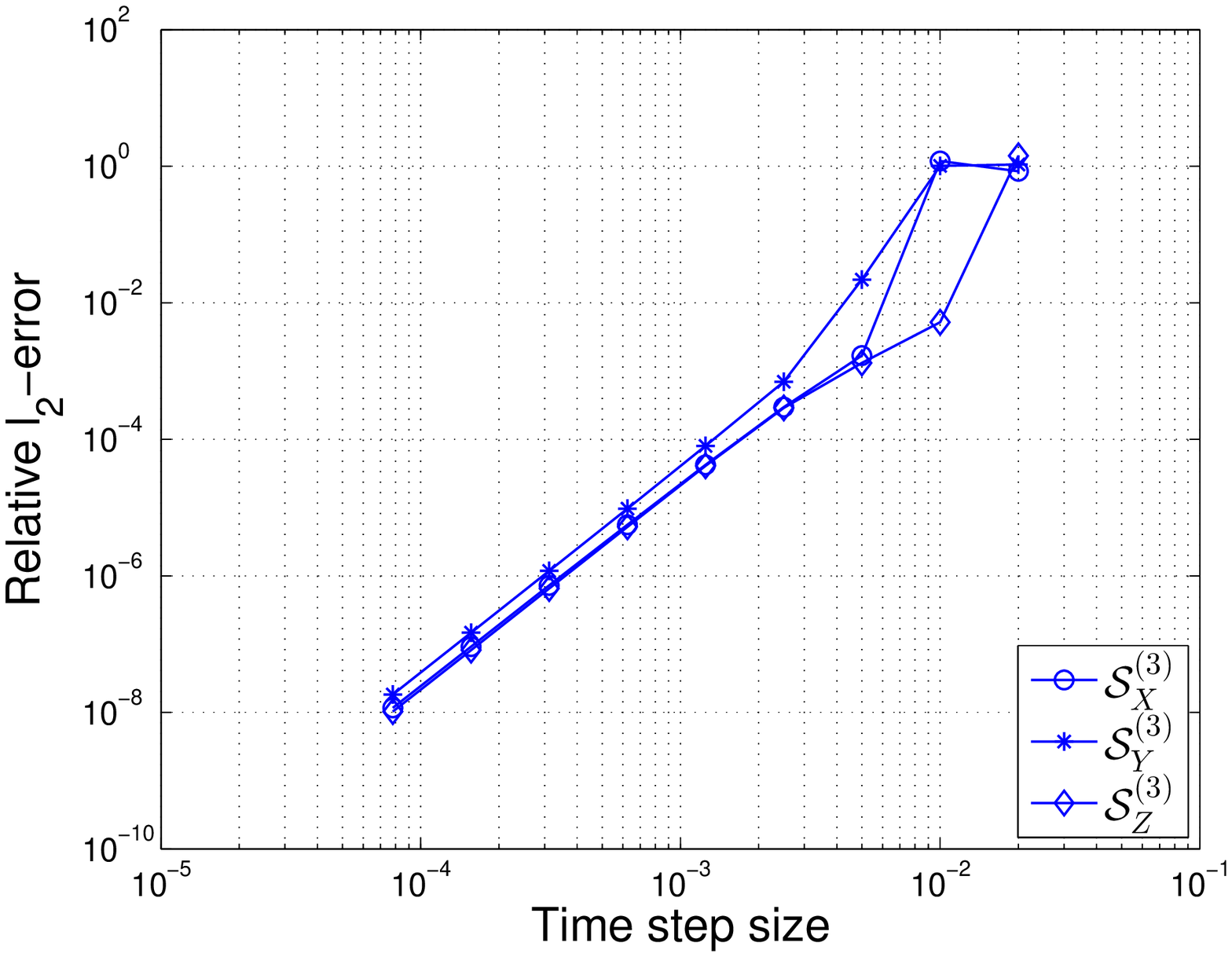} \\
(a)
\end{minipage}
\begin{minipage}{0.5\linewidth}
\centering
\includegraphics[width=2.5in]{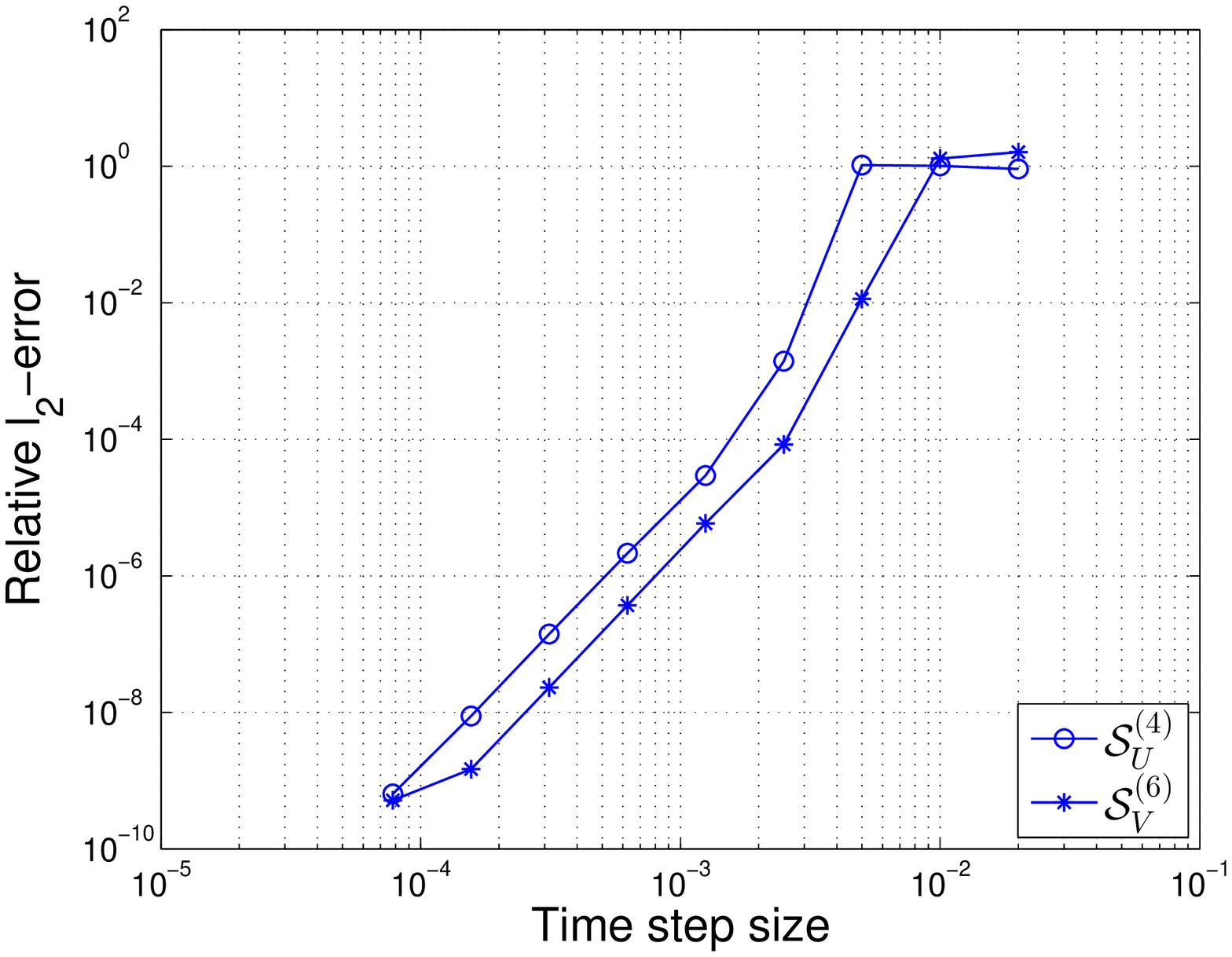} \\
(b)
\end{minipage}
\caption{Relative $l_2$ errors of the traveling wave solution $\phi(x,T_f)$ at $T_f=1/s$ by
(a) the third order methods ${\mathcal S}_X^{(3)}$,  ${\mathcal S}_Y^{(3)}$,
${\mathcal S}_Z^{(3)}$
(b) the fourth order schemes ${\mathcal S}_U^{(4)}$, ${\mathcal S}_V^{(6)}$}
\label{acc_test34}
\end{figure}

We implement the proposed third order schemes ${\mathcal S}_X^{(3)},
{\mathcal S}_Y^{(3)},  {\mathcal S}_Z^{(3)}$ in \eqref{S3XYZ}
and the fourth order schemes ${\mathcal S}_U^{(4)}$ in \eqref{s4w}
and ${\mathcal S}_V^{(6)}$ in \eqref{s6w}.
We set the spectral grid size $h=2^{-5}$, the cut-off limit $K_{tol} = 10^9$
and compare the numerical solutions for various time step sizes $\Delta t$
with the analytic solution~\eqref{wave1d} with $\epsilon=0.03\sqrt{2}$.
Figure~\ref{acc_test34} numerically indicates that
the proposed methods have the third and the fourth order accuracy, respectively.

\subsection{Convergence of the spinodal decomposition problem in 3D}
\begin{figure}[htbp]
\begin{minipage}{0.32\linewidth}
\centering
\includegraphics[width=1.7in]{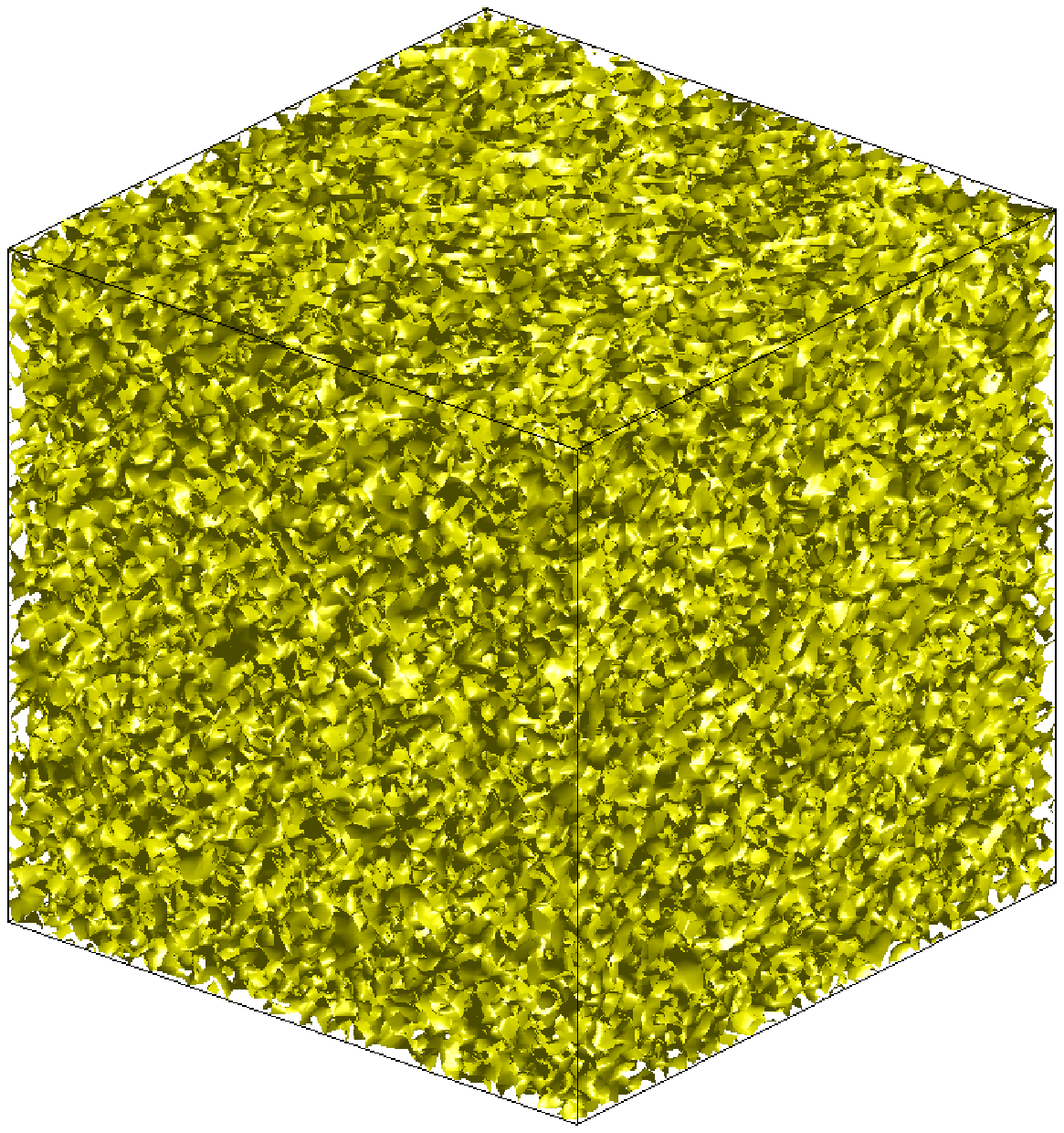} \\
$t=0$
\end{minipage}
\begin{minipage}{0.33\linewidth}
\centering
\includegraphics[width=1.7in]{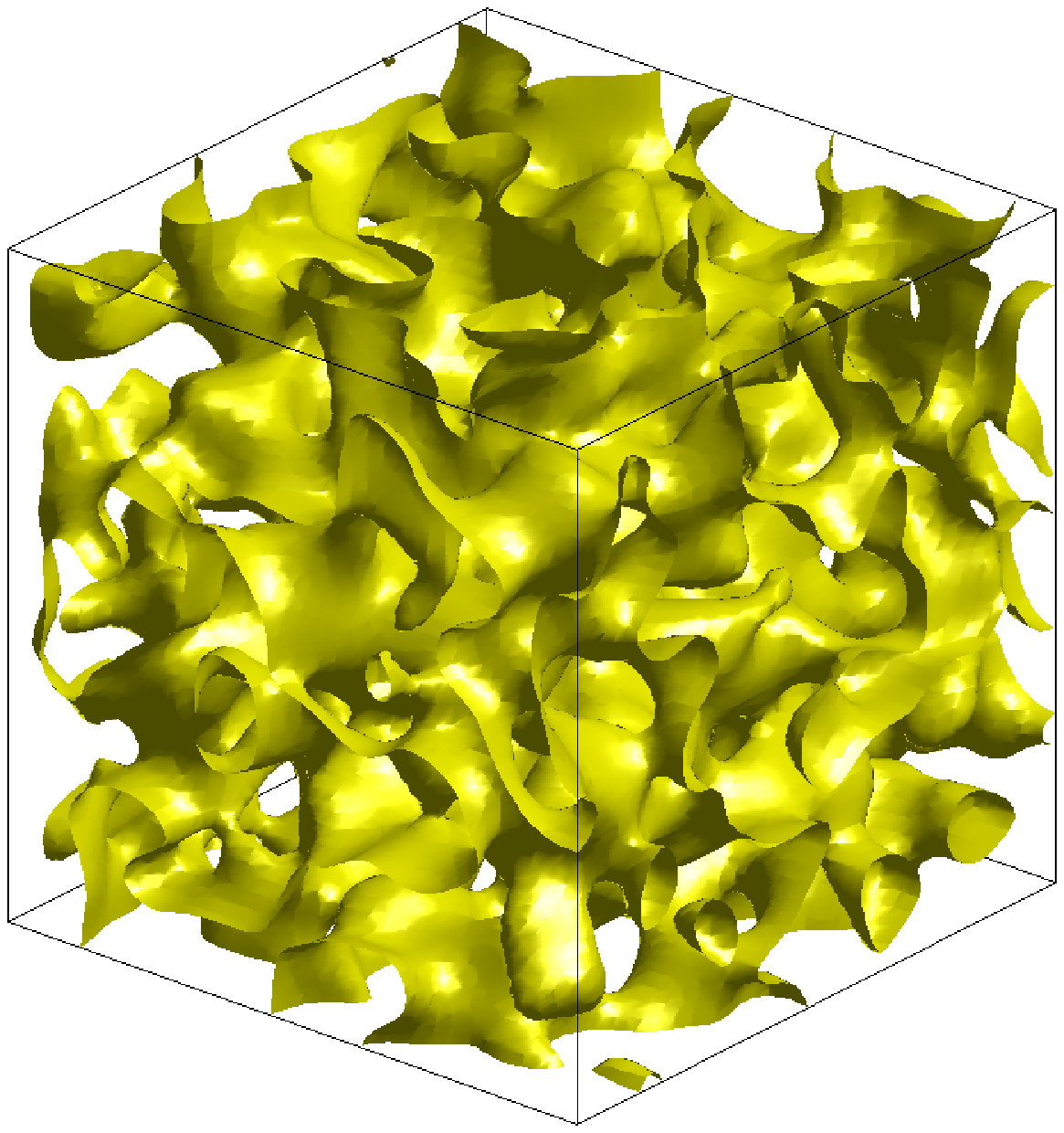} \\
$t=10^{-3}$
\end{minipage}
\begin{minipage}{0.33\linewidth}
\centering
\includegraphics[width=1.7in]{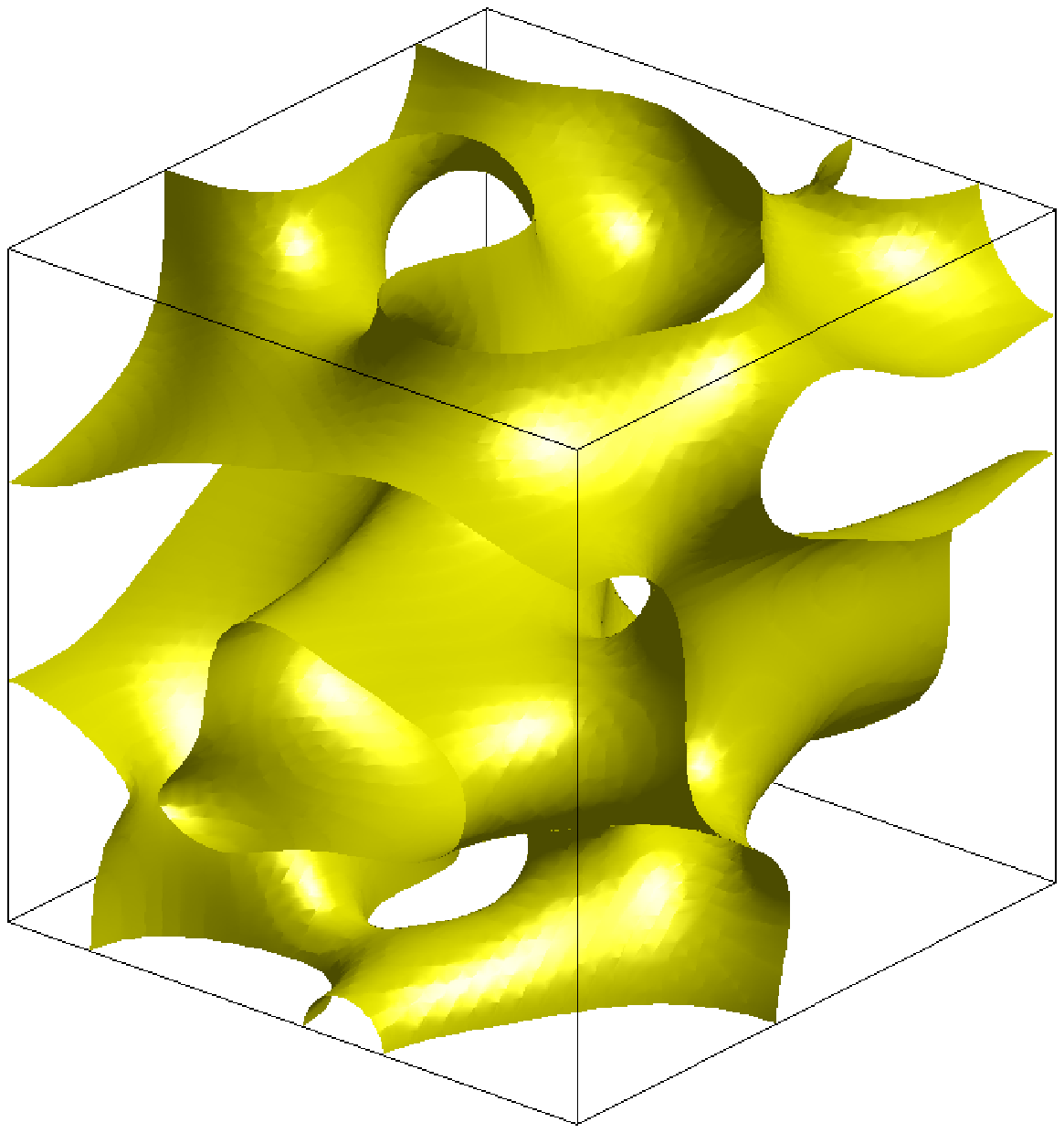} \\
$t=10^{-2}$
\end{minipage}
\caption{The reference solutions $\phi(x,y,z,t)$ by the fourth order
method ${\mathcal S}_V^{(6)}$ with $K_{tol} = 10^9$, and $\Delta
t=10^{-3} / 2^8$.} \label{phase_3d}
\end{figure}

In this subsection, we compute a spinodal decomposition problem
satisfying the AC equation~\eqref{AC_eq} in three-dimensional space with $\epsilon=0.015$.
The intervals $({-}1,{-}1/\sqrt{3})$ and $(1/\sqrt{3},1)$
where $F''(\phi)>0$ are called metastable intervals and
$({-}1/\sqrt{3},1/\sqrt{3})$ where $F''(\phi)<0$ is called the spinodal interval~\cite{Fife}.
It is known that $\phi$ which lies in the spinodal interval is very unstable and
the growth of instabilities results in phase separation,
which is called spinodal decomposition.
In order to check the numerical convergence,
we integrate $\phi(x,y,z,t)$ up to time $T_f=0.01$ by the proposed numerical schemes
with various time step sizes $\Delta t=10^{-3}/2,\cdots,10^{-3}/2^7$.
The initial condition is given on the computational grid
with $h=2^{-6}$ in the domain $\Omega=[0,1]\times[0,1]\times[0,1]$
as $\phi(x,y,z,0) = 0.005\cdot \mbox{rand}(x,y,z)$
where $\mbox{rand}(x,y,z)$ is a random number between $-1$ and $1$.
Figure~\ref{phase_3d} shows the initial and the
reference solutions at $t=10^{-3}, 10^{-2}$
computed by the fourth order numerical scheme ${\mathcal S}_V^{(6)}$
with the numerical parameters $K_{tol} = 10^9$ and $\Delta t=10^{-3} / 2^8$.

\begin{figure}[htbp]
\begin{minipage}{0.49\linewidth}
\centering
\includegraphics[width=2.5in]{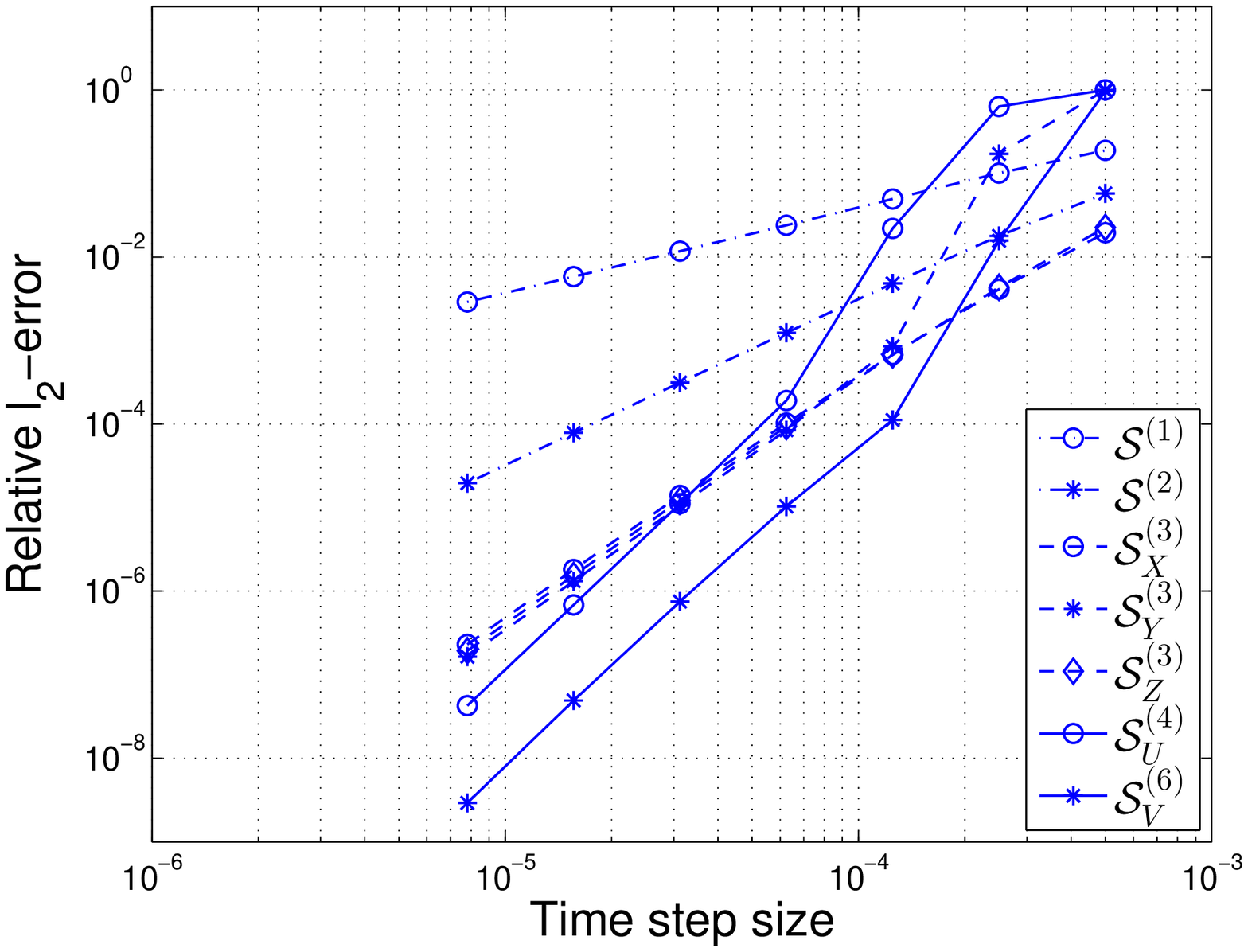} \\
(a)  $K_{tol} = 10^4$
\end{minipage}
\begin{minipage}{0.5\linewidth}
\centering
\includegraphics[width=2.5in]{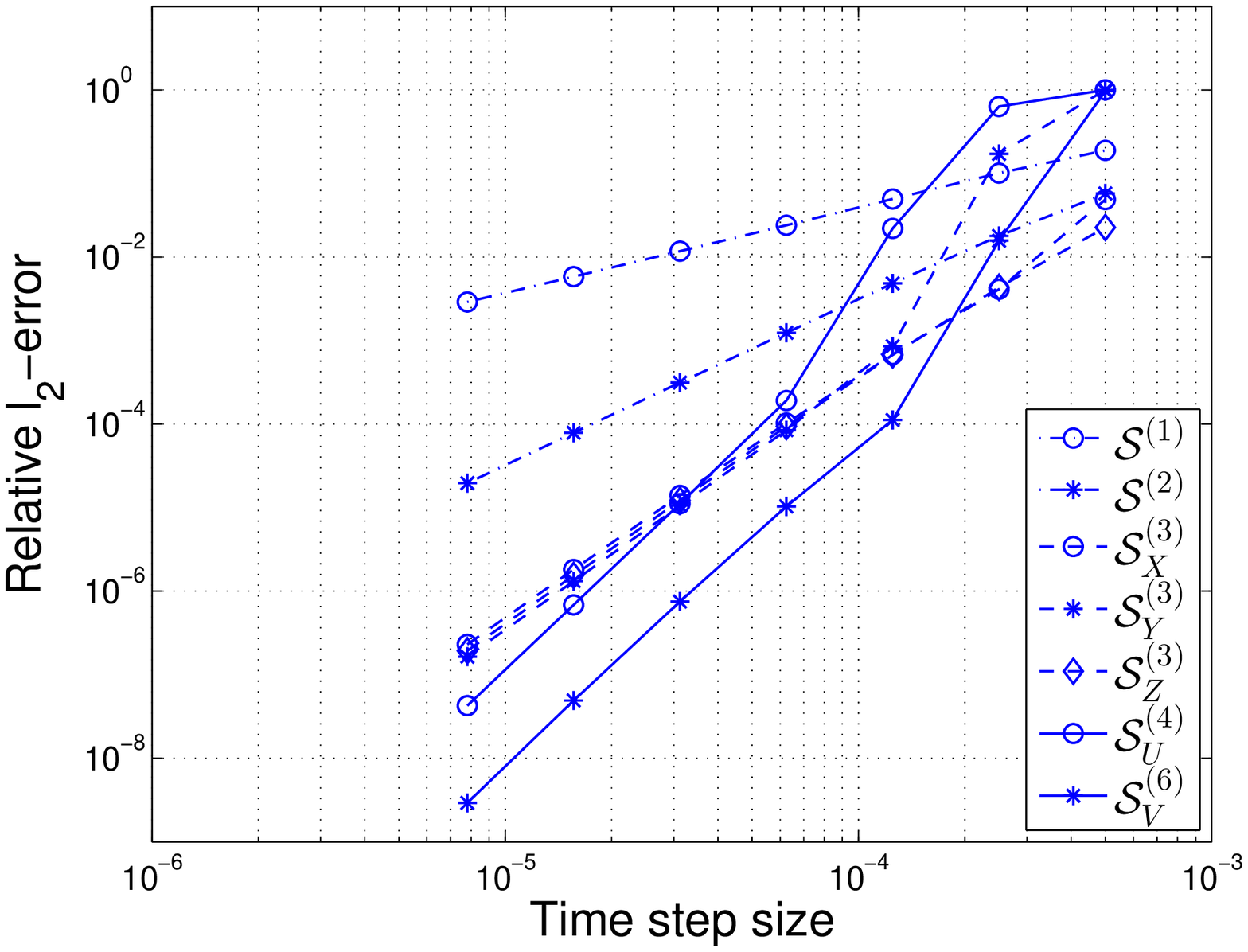} \\
(b) $K_{tol} = 10^9$
\end{minipage}
\caption{Relative $l_2$ errors of $\phi(x,y,z,T_f=0.01)$ by
$\mathcal S^{(1)}$, $\mathcal S_{\omega=1}^{(2)}$, ${\mathcal S}_X^{(3)},
{\mathcal S}_Y^{(3)}, {\mathcal S}_Z^{(3)}$, ${\mathcal
S}_U^{(4)}$, ${\mathcal S}_V^{(6)}$ with
various time step sizes $\Delta t=10^{-3}/2,\cdots,10^{-3}/2^7$.}
\label{acc_test}
\end{figure}

We also implement the first order scheme $\mathcal S^{(1)}$ in \eqref{s1},
the second order scheme $\mathcal S_{\omega=1}^{(2)}$ in \eqref{s2w},
the proposed third order schemes ${\mathcal S}_X^{(3)},
{\mathcal S}_Y^{(3)},  {\mathcal S}_Z^{(3)}$ in \eqref{S3XYZ},
and the fourth order schemes ${\mathcal S}_U^{(4)}$ in \eqref{s4w}
and ${\mathcal S}_V^{(6)}$ in \eqref{s6w}.
The numerical results in Figure~\ref{acc_test} show that
the cut-off value $K_{tol}$ does not play a role when
$\Delta t$ is smaller than $\epsilon^2$ while
the computational results have marginal difference when
$\Delta t$ is greater than $\epsilon^2$.
The accuracy results numerically demonstrate the proposed schemes
provide the expected order of convergence in time.

\section{Conclusions} \label{discu}

We proposed and studied the higher order operator splitting Fourier spectral methods
for solving the AC equation. The methods decompose the AC equation into the
subequations with the heat and the free-energy evolution terms. Unlike the
first and the second order methods, each of the heat and the free-energy evolution
operators has at least one backward evaluation in the higher order methods.
For the third order method, we suggested the three values
$\omega_X,\omega_Y,\omega_Z$ at which $\max \{ |a_j|, |b_j| \}$
have local minimums and we then obtained smaller error than
other $\omega$ values. For the fourth order method, we used two
symmetric combinations of the second order operators.
And a simple cut-off function could limit exponential amplification
of the high frequency modes in the heat operator
and it worked well with the proposed schemes.
We numerically demonstrated, using the traveling wave solution
and the spinodal decomposition problem with random initial values,
that the proposed methods have the third and
the fourth order convergence as expected.

\section*{Acknowledgment}

This research was supported by Basic Science Research Program
through the National Research Foundation of Korea(NRF)
funded by the Ministry of Education(2009-0093827, 2012-002298).


\begin{thebibliography}{00}


\bibitem{AC}
S.M. Allen, J.W. Cahn, A microscopic theory for antiphase boundary
motion and its application to antiphase domain coarsening, Acta
Metall. 27 (1979) 1085--1095.

\bibitem{BCM}
M. Bene$\check{\mbox{s}}$, V. Chalupeck$\acute{\mbox{y}}$, K.
Mikula, Geometrical image segmentation by the Allen--Cahn equation,
Appl. Numer. Math. 51 (2004) 187--205.

\bibitem{DB}
J.A. Dobrosotskaya, A.L. Bertozzi, A wavelet-Laplace variational
technique for image deconvolution and inpainting, IEEE Trans. Image
Process. 17 (2008) 657--663.

\bibitem{ESS}
L.C. Evans, H.M. Soner, P.E. Souganidis, Phase transitions and
generalized motion by mean curvature, Commun. Pur.
Appl. Math. 45 (1992) 1097--1123.

\bibitem{KKR}
M. Katsoulakis, G.T. Kossioris, F. Reitich, Generalized motion by
mean curvature with Neumann conditions and the Allen--Cahn model for
phase transitions, J. Geom. Anal. 5 (1995)
255--279.

\bibitem{FP}
X. Feng, A. Prohl, Numerical analysis of the Allen--Cahn equation
and approximation for mean curvature flows, Numer. Math.
94 (2003) 33--65.

\bibitem{YFL}
X. Yang, J.J. Feng, C. Liu, J. Shen, Numerical simulations of jet
pinching-off and drop formation using an energetic variational
phase-field method, J. Comput. Phys. 218 (2006)
417--428.

\bibitem{Koba}
R. Kobayashi, Modeling and numerical simulations of dendritic
crystal growth, Phys. D 63 (1993) 410--423.

\bibitem{KR}
A. Karma, W.-J. Rappel, Quantitative phase-field modeling of
dendritic growth in two and three dimensions, Phys. Rev. E 57 (1998)
4323--4349.

\bibitem{BWB}
W.J. Boettinger, J.A. Warren, C. Beckermann, A. Karma, Phase-field
simulation of solidification, Annu. Rev. Mater. Res. 32 (2002)
163--194.

\bibitem{Eyre}
D.J. Eyre, An unconditionally stable one-step scheme for gradient
systems,
http://www.math.utah.edu/$\sim$eyre/research/methods/stable.ps.

\bibitem{Yang}
X. Yang, Error analysis of stabilized semi-implicit method of
Allen--Cahn equation, Discrete Cont. Dyn. B 11 (2009)
1057--1070.

\bibitem{SY2}
J. Shen, X. Yang, Numerical approximations of Allen--Cahn and
Cahn--Hilliard equations, Discrete Cont. Dyn. A 28 (2010) 1669--1691.

\bibitem{LLJK}
Y. Li, H.G. Lee, D. Jeong, J. Kim, An unconditionally stable hybrid
numerical method for solving the Allen--Cahn equation, Comput. Math.
Appl. 60 (2010) 1591--1606.

\bibitem{LL}
H.G. Lee, J.-Y. Lee, A semi-analytical Fourier spectral method for
the Allen--Cahn equation, Comput. Math. Appl. 68 (2014) 174--184.

\bibitem{Str}
G. Strang, On the construction and comparison of difference schemes,
SIAM J. Numer. Anal. 5 (1968) 506--517.

\bibitem{GK}
D. Goldman, T.J. Kaper, $N$th-order operator splitting schemes and
nonreversible systems, SIAM J. Numer. Anal. 33 (1996) 349--367.

\bibitem{BC}
S. Blanes, F. Casas, On the necessity of negative coefficients for
operator splitting schemes of order higher than two, Appl. Numer.
Math. 54 (2005) 23--37.

\bibitem{LL2}
H.G. Lee, J.-Y. Lee, A second order operator splitting method
for Allen--Cahn type equations with nonlinear source terms,
Submitted.

\bibitem{LSL}
H.G. Lee, J. Shin, J.-Y. Lee, First and second order operator
splitting methods for the phase field crystal equation, Submitted.

\bibitem{ME}
G.M. Muslu, H.A. Erbay,
Higher-order split-step Fourier schemes for the generalized
nonlinear Schr\"odinger equation,
Math. Comput. Simul. 67 (2005) 581--595.

\bibitem{McL}
R. McLachlan, Symplectic integration of Hamiltonian wave equations,
Numer. Math. 66 (1994) 465--492.

\bibitem{NTK}
N. Ahmed, T. Natarajan, K.R. Rao, Discrete cosine transform, IEEE
Trans. Comput. C-23 (1974) 90--93.

\bibitem{Fife}
P.C. Fife, Models for phase separation and their mathematics,
Electron. J. Diff. Equ. 2000 (2000) 1--26.


%
%
%
%
%
%
%
%
%
%
%
%
%
%

%
%
%



%
%
%
%
%
%


%
%
%

%
%
%

%
%













\end{thebibliography}
\end{document}